\documentclass[a4paper, 12pt, twoside]{article}
\usepackage{amsmath}
\usepackage{amsthm}
\usepackage{amsfonts}
\usepackage{pifont}
\usepackage{bbm}
\usepackage{mathrsfs}
\usepackage{fancyhdr}
\usepackage{graphicx}
\usepackage{color}
\usepackage{psfrag}
\usepackage{time}
\usepackage{txfonts}

\usepackage{stmaryrd}
\usepackage{mathrsfs}
\usepackage{txfonts}
\usepackage{amsfonts}

\usepackage{indentfirst}

\newcommand{\xiaosan}{\fontsize{15pt}{22pt}\selectfont}
\newcommand{\sihao}{\fontsize{14pt}{21pt}\selectfont}

\newcommand{\xiaosi}{\fontsize{12pt}{18pt}\selectfont}

\usepackage{curves}

\numberwithin{equation}{section}

\newtheorem{theorem}{ {Theorem}}[section]

\newtheorem{remark} {   {Remark}}[section]
\newtheorem{corollary} {  {Corollary}}[section]
\newtheorem{lemma} {  {Lemma}}[section]
\newtheorem{definition}{  {Definition}}[section]

\begin{document}
\setlength{\parindent}{2em}
\newpage
\fontsize{12}{22}\selectfont\thispagestyle{empty}
\renewcommand{\headrulewidth}{0pt}
 \lhead{}\chead{}\rhead{} \lfoot{}\cfoot{}\rfoot{}
\noindent

\title{\bf Spectral Gaps of Almost Mathieu Operator in Exponential Regime }
\author{{Wencai  Liu and Xiaoping Yuan*}\\
{\em\small School of Mathematical  Sciences}\\
{\em\small Fudan University}\\
{\em\small  Shanghai 200433, People's Republic of China}\\
{\small 12110180063@fudan.edu.cn}\\
{\small *Corresponding author: xpyuan@fudan.edu.cn}}
\date{}
\maketitle

\renewcommand{\baselinestretch}{1.2}
\large\normalsize
\begin{abstract}
For almost Mathieu operator $(H_{\lambda,\alpha,\theta}u)_n=u_{n+1}+u_{n-1}+2\lambda \cos2\pi(\theta+n\alpha)u_n$,
the dry version of Ten Martini problem predicts that   the spectrum   $\Sigma_{\lambda,\alpha}$ of $ H_{\lambda,\alpha,\theta}$ has all gaps open for all $\lambda\neq 0$
and $ \alpha \in \mathbb{R}\backslash \mathbb{Q}$.
 Avila and  Jitomirskaya  prove   that   $\Sigma_{\lambda,\alpha}$ has all gaps open  for Diophantine $\alpha$  and  $0<|\lambda|<1$.
 In  the present paper, we  show that $\Sigma_{\lambda,\alpha}$ has all gaps open
 for all $ \alpha \in \mathbb{R}\backslash \mathbb{Q}$  with small $\lambda$.
\end{abstract}

\setcounter{page}{1} \pagenumbering{arabic}\topskip -0.82in
\fancyhead[LE]{\footnotesize  Introduction}
\section{\xiaosan \textbf{Introduction  and Main Results  }}
  The almost Mathieu operator (AMO) is the (discrete) quasi-periodic   Schr\"{o}dinger operator  on  $   \ell^2(\mathbb{Z})$:
$$(H_{\lambda,\alpha,\theta}u)_n=u_{n+1}+u_{n-1}+ \lambda v(\theta+n\alpha)u_n,   \text{ with }  v(\theta)=2\cos2\pi \theta,  $$
where $\lambda$ is the coupling, $\alpha $ is the frequency, and $\theta $ is the phase.
\par

 The AMO  is a tight binding model for the Hamiltonian of an electron in a one-dimensional
lattice or
in a two-dimensional lattice, subjecting to a perpendicular (uniform) magnetic field (through a Landau gauge)$ {\cite{Har},\cite{Rau}}$.
This model   also describes a square lattice with anisotropic nearest neighbor coupling and
isotropic next nearest neighbor coupling, or anisotropic coupling to the nearest neighbors and next nearest neighbors
on a triangular lattice$ {\cite{Bel},\cite{Tho}}$.
In addition, the 1980 discovery of the integer Quantum Hall
Effect by von Klitzing$ \cite{Kli} $  leads to a
beautiful theory by Thouless, Kohmoto, Nightingale, and den Nijs.  Central to their theory is the use of the AMO as a model for Bloch electrons in a magnetic field.
   For more applications in physics,  we refer the reader to
 $\cite{Las}$  and the references therein.
 \par
 Besides its application to some fundamental problems in physics, the AMO itself is also
fascinating  because of its remarkable richness of the related spectral theory. In Barry Simon's list of Schr\"{o}dinger operator problems
for the twenty-first century $   \cite{Sim} $, there are three problems   about the AMO. The spectral theory of AMO
has attracted  many   authors, for example,  Avila-Damanik\cite{AD}, Avila-Jitomirskaya\cite{AJ1}, \cite{AJ2},  Avron-Simon \cite{AS},\cite{Avr}, Bourgain\cite{B2},  Jitomirskaya-Last\cite{JL1} and so on.
\par
Here we are concerned with the topological structure of the spectrum, which is heavily related to the  arithmetic properties of  frequency $\alpha $.  If $\alpha=p/q$ is rational, it is well known that the spectrum consists of the union of $q$ intervals called bands, possibly touching the endpoints. When $\alpha\in\mathbb R\setminus\mathbb Q$ and $\lambda\neq 0$, the spectrum set $\Sigma_{\lambda,\alpha}$ of $H_{\lambda,\alpha,\theta}$ (  in this case the spectrum of $H_{\lambda,\alpha,\theta}$ is independent of $\theta$) has been conjectured for a long time to be a Cantor set. This conjecture is named after the Ten  Martini Problem\footnote{  Ten Martini Problem is  the Fourth Problem   in $\cite{Sim}.$}.    It has been solved by  Avila-Jitomirskaya completely \cite{AJ1}
 by  Anderson localization (i.e., only pure point spectrum with exponentially decaying eigenfunctions) of $H_{\lambda,\alpha,\theta}$
 when $  |\lambda| >e^{ \frac{16\beta}{9}}$,  where
 \begin{equation}\label{G11}
  \beta= \beta(\alpha)=\limsup_{n\rightarrow\infty}\frac{\ln q_{n+1}}{q_n},
 \end{equation} and $ \frac{p_n}{q_n} $ be the continued fraction approximants ($\S 2.5$) to $\alpha$.    See $\cite{AJ1}$  for more historic backgrounds about the Ten  Martini Problem.
Recently, the condition $  |\lambda| >e^{ \frac{16\beta}{9}}$ has been refined to $  |\lambda| >e^{ \frac{3\beta}{2}}$ by the present authors  $\cite{LIU}$.

\par
About the  topological structure of the spectrum of the AMO, a stronger conjecture is so-called the dry version of the Ten Martini Problem. In order to state it, we introduce the integrated density of states $N_{\lambda,\alpha}(E)$ (see $(\ref{xin})$) of AMO, which
 is a continuous non-decreasing surjective function with
$N_{\lambda,\alpha}\text{:\;}\mathbb{R} \mapsto[0,1]$.  For  $ \alpha\in \mathbb{R}\backslash\mathbb{ Q}$, the basic relation between $\Sigma_{\lambda,\alpha}$
%\footnote{For  $\alpha\in \mathbb{R} \backslash
%\mathbb{Q}$,  denote $ \Sigma_{\lambda,\alpha} $   by  the spectrum  of $H_{\lambda,\alpha,\theta}$ since   $\Sigma_{\lambda,\alpha, \theta }$ does not depend on $\theta$ ($\S 2.3$), where $ \Sigma_{\lambda,\alpha, \theta}$ is the spectrum of $H_{\lambda,\alpha,\theta}$. For   $\alpha \in \mathbb{Q}$, let
 %$ \Sigma_{\lambda,\alpha}=\bigcup_{\theta} \Sigma_{\lambda,\alpha, \theta }  $.  }
 and  $N_{\lambda,\alpha}$
is that $E\notin \Sigma_{\lambda,\alpha}$ if and only if $N_{\lambda,\alpha}$ is constant in a neighborhood of $E$.
Each   connected component of  $\mathbb{R}\backslash \Sigma_{\lambda,\alpha}$  is called a    gap of  $\Sigma_{\lambda,\alpha}$.
If $E$ is an endpoint of some   gap, then $N_{\lambda,\alpha}(E)\in \alpha \mathbb{Z}+\mathbb{Z}$ (combining $ \cite{Her} $ with $ \cite{Jo} $).
The dry version of  Ten Martini Problem predicts the  converse is also true.  Concretely, $N_{\lambda,\alpha}(E)\in \alpha \mathbb{Z}+\mathbb{Z}$   with $E\in \Sigma_{\lambda,\alpha}$ implies   $E$ is an endpoint
of some   gap
   for all $\lambda\neq 0$ and  $\alpha\in\mathbb{R}\backslash \mathbb{Q}$ (this obviously implies the Ten Martini Problem).
   For convenience,
we say that all gaps of $\Sigma_{\lambda,\alpha}$ are open if  $N_{\lambda,\alpha}(E)\in \alpha \mathbb{Z}+\mathbb{Z}$   and $E\in \Sigma_{\lambda,\alpha}$ implies   $E$ is an endpoint
of some   gap. Equivalently, {\it the dry version of the Ten Martini Problem  predicts $\Sigma_{\lambda,\alpha}$ has all gaps open for
all $\lambda\neq 0$ and $\alpha\in\mathbb{R}\backslash \mathbb{Q}$.}
\par

In proving the dry version of the   Ten Martini Problem, much progress has been recently achieved by many authors. The proofs depend on whether $\beta(\alpha)>0$ or $\beta(\alpha)=0$.  One usually calls set $\{\alpha \in \mathbb{R}\backslash \mathbb{Q}|\;\beta(\alpha)>0\}$    exponential regime and set $\{\alpha\in \mathbb{R}\backslash \mathbb{Q}|\;\beta(\alpha)=0\}$   sub-exponential regime.

 \par

In the exponential regime ($\beta(\alpha)>0$), for any  $   \varepsilon>0$, one has
$$|\frac{p_n}{q_n}-\alpha|\leq\frac{1}{q_nq_{n+1}}\leq \frac{1}{e^{( \beta -\varepsilon )q_n}}$$ if
  $n$ is large enough (see $(\ref{XG29})$). This means rational number $\frac{p_n}{q_n}$  is exponentially close to $\alpha$, thus
the   gaps of spectrum  $\Sigma_{\lambda,\alpha}$  can be   rational approximated  by the gaps  of $\Sigma_{\lambda, \frac{p_n}{q_n}}$\footnote{For   $\alpha \in \mathbb{Q}$,
 $ \Sigma_{\lambda,\alpha}=\bigcup_{\theta} \Sigma_{\lambda,\alpha, \theta }  $, where $ \Sigma_{\lambda,\alpha, \theta}$ is the spectrum of $H_{\lambda,\alpha,\theta}$.  }.
 Choi, Elliott, and Yui  in $ \cite{CEY}$ set up   $1/3$-H\"{o}lder continuity of the spectrum for $\lambda=1$, i.e.,
 $\text{Dist}(\Sigma_{1,\alpha_1},\Sigma_{1,\alpha_2})<C|\alpha_1-\alpha_2|^{1/3}$, and give a good  estimate for  the gaps of $ \Sigma_{1,\frac{p_n}{q_n}}$, where  $\text{Dist}(K_1,K_2)$ means
Hausdorff distance  between  two subsets $K_1\subset\mathbb R$ and $K_2\subset\mathbb R$.  Then they prove: if $ \lambda=1$,  $\Sigma_{\lambda,\alpha}$ has all gaps open for $\beta>9\ln2+3\ln3$. In $ \cite{Avms}  $, Avron,   Mouche  and  Simon set up
 $1/2$-H\"{o}lder continuity of the spectrum. Combining  $ \cite{Avms}  $ with $ \cite{CEY}$,   Avila and Jitomirskaya   $\cite{AJ1}$
 obtain that  $\Sigma_{\lambda,\alpha}$ has all gaps open for   $ 0<\beta \leq \infty$ and $e^{-\beta} <|\lambda|<e^{\beta}$.
 In particular, for  $\alpha\in \mathbb{R}\backslash \mathbb{Q}$ such that  $\beta =\infty$, $\Sigma_{\lambda,\alpha}$ has all gaps open if  $\lambda\neq 0$. Thus in the present paper,
   unless stated otherwise,   we always  assume $\alpha\in \mathbb{R}\backslash \mathbb{Q}$ such that $\beta=\beta(\alpha)<\infty$ and $\lambda\neq 0$.
 \par

\par
Now let us return to the sub-exponential regime ($\beta(\alpha)=0$). Of this case is the well-known Diophantine number.   We say $\alpha \in \mathbb{R}\backslash \mathbb{Q} $ satisfies a Diophantine condition $\text{DC}(\kappa,\tau)$ with $\kappa>0$ and $\tau>0$,
if
$$ |q\alpha-p|>\kappa |q|^{-\tau}  \text{ for   any } (p,q)\in \mathbb{Z}^2, q\neq 0.$$
Let $\text{DC}=\cup_{\kappa>0,\tau>0}\text{DC}( \kappa,\tau)$. We say $\alpha$  satisfies   Diophantine condition, if $\alpha\in \text{DC}$.
 Notice that the set DC is a real subset of the sub-exponential regime, i.e., $\text{DC} \subsetneqq\{\alpha:\beta(\alpha)=0\}$.
For   $ \alpha \in DC$, Puig in $ \cite{P1}$ and $ \cite{P2}$ develops a way to estimate   the gaps via establishing   reducibility ($\S 2.1$).
He proves: for  $ \alpha \in DC$, if $N_{\lambda,\alpha}(E)\in \alpha \mathbb{Z}+\mathbb{Z}$ with $E\in \Sigma_{\lambda,\alpha}$ and cocycle $S_{\lambda,E}$ ($\S 2.3$)
 is analytically reducible, then $E$ is an endpoint of some   gap.
  In $ \cite{AJ2}$, Avila and Jitomirskaya develop a quantitative version of Aubry duality  and use it to obtain a sharp estimate  of the rotation
 number $\rho(\alpha,A)$ ($\S 2.2$) with $A=S_{\lambda,E}$ for   $\alpha\in DC$. As a result, they  establish   reducibility for the cocycle $S_{\lambda,E}$ when $0<|\lambda|<1$ with $E \in \Sigma_{\lambda,\alpha}$ and       $N_{\lambda,\alpha}(E)\in \alpha \mathbb{Z}+\mathbb{Z}$.
  Combining with Puig's discussion,  they      show that  $\Sigma_{\lambda,\alpha}$ has all gaps open if
 $ \alpha \in DC$ and $ 0<|\lambda|<1$.
 \par
 In conclusion, we give a list for the unsolved cases about the dry version of the Ten Martini Problem\footnote{By Aubry duality, it suffices to discuss $0<|\lambda|\le 1$.}:
 \begin{enumerate}

  \item $\quad \alpha\in\{\alpha:\;0<\beta(\alpha)<\infty\},\; 0<|\lambda|\le e^{-\beta}$;
 \item $\quad \alpha\in\text{DC},\; |\lambda|=1$;
 \item $\quad \alpha\in\{\alpha:\;\beta(\alpha)=0\}\setminus DC,\; 0<|\lambda|\le 1$.

 \end{enumerate}
 \par
 In the present paper, we prove the following theorem:
\begin{theorem} (\textbf{Main Theorem})\label{Th11}
For every $ \alpha \in \mathbb{ R}\backslash \mathbb{Q}$ such  that $\beta=\beta(\alpha)<\infty$, there exists  a absolute constant $C $, such that  $\Sigma_{\lambda,\alpha}$ has all gaps open if
  $0 <|\lambda|< e^{-C \beta}$.
 \end{theorem}

  \begin{remark} The main contribution in the present paper is that the unsolved regime $|\lambda|\in(0,e^{-\beta}]$ in  case 1 is shrunk to $|\lambda|\in[e^{-C\beta}, e^{-\beta}]$.  We should point out that the constant $C>0$ is very large. Therefore there is a long way to decrease it to $C<1$ such that the problem is solved completely.  The unsolved case 3 is now solved  by letting $\beta=0$ in Theorem \ref{Th11}, except $|\lambda|=1$. Actually, the case 3  is  solved by careful checking the proofs of  \cite{AJ2} and \cite{P2}.
    \end{remark}

  The present paper is organized as follows:

  \par
  In \S2, we give some preliminary notions and facts which are taken from \cite{AJ1}.
  \par
 In \S3, we obtain the strong localization estimate    of the Aubry dual model $\hat H_{\lambda  ,\alpha,\theta} $
 in the exponential regime (i.e., $\beta(\alpha)>0$).
 \par
In \S4,     we  set up    sharp   estimate of
 the rotation number  $(\text{Theorem } \ref{Th415})$ for resonant phase by developing   the  quantitative version of Aubry duality  in exponential regime. This process is the same as to set up almost reducibility for  cocycles $S_{\lambda,E}$.
\par
 In \S5,
 we  obtain the  analytic reducibility in a trip domain for non-resonant phase ( Theorem \ref{Th417}) by   constructing a new
     reducible matrix  in  $\text{PSL}(2,\mathbb{R})$ (by  Lemma  \ref{matrix B} and   Theorem \ref{XTh54}).
   Combining with the sharp   estimate of
   rotation number in \S 4, we    set up the reducibility for   cocycle $S_{\lambda,E}$ when
    $E \in \Sigma_{\lambda,\alpha}$ and $E$ satisfies  $N_{\lambda,\alpha}(E)\in \alpha \mathbb{Z}+\mathbb{Z}$  $(\text{Theorem } \ref{Th420})$.
 \par
 In \S6, in order to use Puig's method, we  generalize his result  to  exponential regime by  KAM iteration  $\;(\text{Theorem } \ref{Th51})$.
  In the end, we   give a summary  about the dry version of  Ten Martini Problem (Theorem \ref{Conclution} ).

\section{\xiaosan \textbf{ Preliminaries }}
\subsection{\xiaosan \textbf{  Cocycles}}
Denote by $ \text{SL}(2,\mathbb{C})$ the all complex  $2\times 2 $-matrixes  with determinant 1.
We say a function
  $f\in C^{\omega}(\mathbb{R}/\mathbb{Z},  \mathbb{C})$ if $f$ is well  defined in $ \mathbb{R}/\mathbb{Z}$, i.e., $f(x+1)=f(x)$
  and $f$ is   analytic in
    a neighbor of $\Im x=0$. The definitions of $ \text{SL}(2,\mathbb{R})$ and  $  C^{\omega}(\mathbb{R}/\mathbb{Z},  \mathbb{R})$
    are similar to those of $ \text{SL}(2,\mathbb{C})$ and $C^{\omega}(\mathbb{R}/\mathbb{Z},\mathbb{C})$, respectively, except that the involved matrixes are real and the functions are real   analytic.
A $C^{\omega}$-cocycle in $ \text{SL}(2,\mathbb{C})$  is a pair $(\alpha,A)\in \mathbb{R} \times  C ^{\omega} (\mathbb{R}/\mathbb{Z}, \text{SL}(2,\mathbb{C})) $, where
$A  \in    C ^{\omega} (\mathbb{R}/\mathbb{Z}, \text{SL}(2,\mathbb{C}))$ means
      $A(x)\in \text{SL}(2,\mathbb{C})$      and the elements of $A$ are in $  C^{\omega}(\mathbb{R}/\mathbb{Z},  \mathbb{C})$.
       Sometimes, we say
     $ A$ a $C^{\omega}$-cocycle for short, if there is no ambiguity. Note that all functions, cocycles in the present paper are analytic.
     Thus we often   do not mention the analyticity, for instance, we say $A$  a cocycle instead of   $C^{\omega}$-cocycle.

    The Lyapunov exponent for the cocycle $A$ is given  by
 \begin{equation}\label{G21}
    L(\alpha,A)=\lim_{n\rightarrow\infty} \frac{1}{n}\int_{\mathbb{R}/\mathbb{Z}} \ln \| A_n(x)\|dx,
 \end{equation}
where
\begin{equation}\label{G22}
     A_n(x) = A(x+(n-1)\alpha)A(x+(n-2)\alpha)\cdots A(x).
\end{equation}
By Corollary 2   in   \cite{FUR}  (since irrational rotations are uniquely ergodic)
\begin{equation}\label{G23}
   L(\alpha,A)=\lim_{n\rightarrow\infty} \sup_{x\in\mathbb{R}/ \mathbb{Z}}\frac{1}{n} \ln \| A_n(x)\|,
\end{equation}
that is,  the convergence in ($\ref{G23}$) is  uniform   with respect to  $x\in\mathbb{R}$.
 Precisely, $ \forall \varepsilon >0$,
\begin{equation}\label{G24}
  \| A_n(x)\|\leq e^{(L(\alpha,A)+\varepsilon)n},  \text {for  }n  \text { large enough}.
\end{equation}

Given two cocycles $(\alpha,A)$  and $(\alpha,A^{\prime})$, a conjugacy between them is a   cocycle
$ B \in C ^{\omega}(\mathbb{R}/\mathbb{Z},  \text{SL}(2,\mathbb{C}))$ such that
 \begin{equation}\label{G25}
   B(x+\alpha)^{-1}A(x)B(x)=A^{\prime}.
 \end{equation}
  The notion
of real conjugacy (between real cocycles) is the same as before, except that we ask for
 $ B \in C ^{\omega}(\mathbb{R}/\mathbb{Z}, \text{PSL}(2,\mathbb{R}))$, i.e., $B(x+1)= \pm B(x)$ and $\det B=1$.
We say that cocycle $(\alpha,A)$  is reducible if it is  conjugate to a constant cocycle.

 \subsection{\xiaosan \textbf{The rotation number}}
 Let
$A(\theta)=\left(
             \begin{array}{cc}
               a(\theta) & b(\theta) \\
               c(\theta)& d (\theta)\\
             \end{array}
           \right)
$, we define  the map $T_{\alpha,A}:(\theta,\varphi)\in \mathbb{T}\times \frac{1}{2}\mathbb{T} \mapsto
(\theta+\alpha,\varphi_{\alpha,A}(\theta,\varphi))\in \mathbb{T}\times \frac{1}{2}\mathbb{T},$
 with $\varphi_{\alpha,A}= \frac{1}{2\pi}\arctan (\frac{c(\theta)+d(\theta)\tan2\pi\varphi}{a(\theta)+b(\theta)\tan2\pi\varphi})$,
 where $\mathbb{T}=\mathbb{R}/\mathbb{ Z}$.
Assume now that     A :$ \mathbb{R}/   \mathbb{ Z}\rightarrow \text{SL}(2,\mathbb{R})$ is homotopic to the identity,
then  $T_{\alpha,A} $ admits a continuous lift $\tilde{T}_{\alpha,A}:(\theta,\varphi)\in \mathbb{R}\times  \mathbb{R} \mapsto
(\theta+\alpha,\tilde{\varphi}_{\alpha,A}(\theta,\varphi))\in \mathbb{R}\times  \mathbb{R} $  such that
$\tilde{\varphi}_{\alpha,A}(\theta,\varphi)  \mod  \frac{1}{2}\mathbb{Z}= \varphi _{\alpha,A}(\theta,\varphi  ) $
and $\tilde{\varphi}_{\alpha,A}(\theta,\varphi) -\varphi$ is well defined on $\mathbb{T}\times \frac{1}{2}\mathbb{T}$.
 The number $\rho(\alpha,A)=\limsup_{n\rightarrow \infty}\frac{1}{n}(p_2\circ \tilde{T}^n_{\alpha,A}(\theta,\varphi)-\varphi)\mod \frac{1}{2}\mathbb{Z},$
 does not depend on the choices of $ \theta$ and $\varphi$,
where $p_2(\theta,\varphi)=\varphi,$ and is called the   rotation number
of $(\alpha,A)$ $ {\cite{Her} ,\cite{JM}}$.
\par
It  follows  from the definition that ( p.8, \cite{AJ2})
\begin{equation}\label{XXG26}
  ||\rho(\alpha,A)-\theta||_{\mathbb{R}/2\mathbb{Z}}<C\sup_{x\in \mathbb{R}}|| A(x)-R_\theta||,
\end{equation}
where $||x||_{\mathbb{R}/2\mathbb{Z}}=\min_{\ell\in \mathbb{Z}}|x-\frac{\ell}{2}|$ and $||\cdot||$ is any Euclidean norm,
and
$$
R_\theta=
\left(
       \begin{array}{cc}
         \cos2\pi \theta & -\sin2\pi \theta \\
          \sin2\pi \theta &  \cos2\pi \theta \\
       \end{array}
     \right).
$$
\par
If $A, A':\mathbb{R} / \mathbb{ Z}\mapsto \text{SL}(2,\mathbb{R})$   and
$ B :\mathbb{R}/ 2\mathbb{Z}\mapsto \text{SL}(2,\mathbb{R})$  (notice that $ B :\mathbb{R}/  \mathbb{Z}\mapsto \text{PSL}(2,\mathbb{R})$  implies  $ B :\mathbb{R}/ 2 \mathbb{Z}\mapsto \text{SL}(2,\mathbb{R})$ )   such
that $A $ is  homotopic to the identity and $ B(x+\alpha)^{-1}A(x)B(x) =A^{\prime}$, then   $ A^\prime$ is  homotopic to the identity
and
$2\rho(\alpha,A)-2\rho(\alpha,A^\prime)=  k\alpha \mod  \mathbb{Z } $,  where  $k$  is the degree of $B$ ( denoted by $\deg (B)$), i.e., $x\mapsto B(x)$ is homotopic to
$x\mapsto R_{\frac{kx}{2}} $.

 \subsection{\xiaosan \textbf{Almost Mathieu  cocycles and  the integrated density of states }}
For the almost Mathieu operators  \{$ H_{\lambda,\alpha,\theta} \}_{\theta\in \mathbb{R}}$,  the spectrum  of operator
 $ H_{\lambda,\alpha,\theta}  $
 does not depend on $ \theta$,  denoted by  $  \Sigma_{\lambda,\alpha}$. Indeed, shift is an unitary
 operator on $\ell^2(\mathbb{Z})$, thus $  \Sigma  _{\lambda,\alpha,\theta}= \Sigma  _{\lambda,\alpha,\theta+\alpha}$, where $  \Sigma _{\lambda,\alpha,\theta}$
 is the spectrum of  $ H_{\lambda,\alpha,\theta}  $. By  the minimality of $ \theta\mapsto \theta+\alpha$ and  continuity of  spectrum
  $  \Sigma  _{\lambda,\alpha,\theta}$ with respect  to  $\theta$, the statement follows.
  \par
 Let
 $$
 S_{\lambda,E}=
 \left(
     \begin{array}{cc}
       E-2\lambda \cos2\pi x & -1 \\
       1 & 0 \\
     \end{array}
   \right).
   $$
We call $(\alpha,S_{\lambda,E})$ almost Mathieu cocycle.   It's easy to see that almost Mathieu cocycle is  homotopic to the identity,
  and let
$\rho_{\lambda,\alpha}(E)\in[0,\frac{1}{2}]$ be the rotation number of the almost Mathieu cocycle $(\alpha,S_{\lambda,E})$.
\par
Next we will give the definition of the integrated density of states $N_{\lambda,\alpha}$, which has been   mentioned in \S 1.
\par
Let $H$ be a bounded self-adjoint operator on $\ell ^2(\mathbb{Z})$. Then
  $(H-z)^{-1}$ is analytic
in  $ \mathbb{C} \backslash \Sigma( H)$, where $\Sigma( H)$ is the spectrum of $H$,
and we have for $ f\in \ell^2$
\begin{equation*}
  \Im\langle(H-z)^{-1}f,f\rangle=\Im z \cdot ||(H-z)^{-1}f||^2,
\end{equation*}
where $\langle\cdot,\cdot\rangle$ is the usual inner product in $\ell ^2(\mathbb{Z})$.
Thus
\begin{equation*}
 \phi_f(z)=   \langle(H-z)^{-1}f,f\rangle
\end{equation*}
is an analytic function on the upper half plane with $\Im \phi_f\geq 0$ ( $\phi_f$ is a so-called Herglotz function).
\par
Therefore one has a representation
\begin{equation}
 \phi_f(z)=   \langle(H-z)^{-1}f,f\rangle=\int_{\mathbb{R}} \frac{1}{x-z}d\mu^{f}(x)
\end{equation}
where $\mu^{f} $ is the spectral measure associated to $f$.
\par
Fix almost  Mathieu  operator $  H_{\lambda  ,\alpha,\theta}   $. Denote by $\mu^{f}_{\lambda  ,\alpha,\theta} $
the spectral measure of operator  $  H_{\lambda  ,\alpha,\theta}   $ and vector $f$ as before.
The  integrated density of states (IDS) $N_{\lambda  ,\alpha}$ is obtained by   averaging the spectral measure  $\mu_{\lambda , \alpha,\theta}^{\delta_0}$ with respect to $\theta$, i.e.,
\begin{equation} \label{xin}
  N_{\lambda ,\alpha}(E)=\int_{\mathbb{R}/\mathbb{Z}} \mu^{\delta_0}_{\lambda  ,\alpha,\theta}(-\infty,E]d\theta,
\end{equation}
where $\delta_0$ is the normal vector in $\ell^2(\mathbb{Z})$ with 0th  component being 1, others being 0.
\par
Between the integrated density  of   states $N_{\lambda,\alpha}(E)$ and the rotation number  $\rho_{\lambda,\alpha}(E)$,
there is   the following relation $ {\cite{Jo1}}$:
 \begin{equation}\label{G26}
    N_{\lambda,\alpha}(E)=1-2\rho_{\lambda,\alpha}(E).
 \end{equation}
 In particular, $  N_{\lambda,\alpha}(E)\in \alpha \mathbb{Z}+\mathbb{\mathbb{Z}}$ is equivalent to $2\rho_{\lambda,\alpha}(E)\in \alpha \mathbb{Z}+\mathbb{\mathbb{Z}}$.
 \par
Let $L_{\lambda,\alpha}(E)=L(\alpha,S_{\lambda,E})$ be the Lyapunov exponent of $S_{\lambda,E}$.
In $\cite{BJ1} $ Bourgain and Jitomirskaya obtain the accurate value of Lyapunov exponent when $E\in \Sigma_{\lambda,\alpha}$.
 \begin{theorem} $ (\cite{BJ1})$ \label{Th21}
 For every $\alpha\in\mathbb{R} \backslash\mathbb{Q}$, $\lambda\in \mathbb{R}$ and $E\in \Sigma_{\lambda,\alpha}$,  one has
 $L_{\lambda,\alpha}(E)=\max\{\ln|\lambda|,0\} $.

\end{theorem}

 \subsection{ Classical Aubry duality}

 Let $ \hat{H} _{\lambda,\alpha,\theta}   =  \lambda H_{\lambda^{-1},\alpha,\theta}   $. If $\alpha\in\mathbb{R} \backslash  \mathbb{Q}$,  then the spectrum of $ \hat{H}_{\lambda,\alpha,\theta} $     is exactly $\Sigma_{\lambda,\alpha}$  \cite{GJLS}. $  \hat{H} _{\lambda,\alpha,\theta}$
 is called   Aubry dual model of $ H_{\lambda ,\alpha,\theta} $.
Classical Aubry duality expresses an algebraic relation between the families of
operators $ \{\hat{H} _{\lambda,\alpha,\theta}\}_{\theta\in\mathbb{R}} $  and $ \{ {H} _{\lambda,\alpha,x}\}_{x\in\mathbb{R}} $
by  Bloch waves, i.e., if
 $ u:\mathbb{R}/\mathbb{Z}\mapsto\mathbb{C}$ is an $L^2$  function whose Fourier coefficients  $\hat u$  satisfy
  $ \hat{H} _{\lambda,\alpha,\theta}\hat{u}=E\hat{u}$, then

  $$U(x)= \left(
           \begin{array}{c }
             e^{2\pi i \theta }u(x) \\
             u(x-\alpha)\\
           \end{array}
         \right)
  $$
satisfies $S_{\lambda,E}(x)\cdot U(x)=e^{2\pi i \theta}U(x+\alpha).$
\subsection{Continued fraction expansion}
Define as usual for $0\leq\alpha<1,$
$$ a_0=0,\alpha_0=\alpha,$$
and inductively for $k>0,$
$$a_k=\lfloor \alpha_{k-1}^{-1}\rfloor, \alpha_k=\alpha_{k-1}^{-1}-a_k,$$
where $\lfloor t \rfloor$ denotes the greatest integer less than or equal $t$.
\par
We define
$$
\begin{array}{cc}
            p_0=0, & q_0=1, \\
            p_1=1,& q_1=a_1 ,
          \end{array}
$$
and inductively,
\begin{eqnarray*}
  p_k &=& a_k p_{k-1}+p_{k-2}, \\
  q_k &=& a_k q_{k-1}+q_{k-2}.
\end{eqnarray*}
Recall that     $\{q_n\}_{n\in \mathbb{N}}$ is the sequence of best denominators of irrational number $\alpha$,
since it satisifies
\begin{equation}\label{XG28}
\forall 1\leq k <q_{n+1}, \| k\alpha\|_{\mathbb{R}/\mathbb{Z}}\geq ||q_n\alpha||_{\mathbb{R}/\mathbb{Z}},
\end{equation}
where $||x||_{\mathbb{R}/\mathbb{Z}}=\min_{\ell\in \mathbb{Z}}|x-\ell|$. Moreover, we also have the following estimate,

\begin{equation}\label{XG29}
      \frac{1}{2q_{n+1}}\leq\Delta_n\triangleq \|q_n\alpha\|_{\mathbb{R}/\mathbb{Z}}\leq\frac{1}{q_{n+1}}.
\end{equation}

\section{ Strong localization estimate for $0<\beta(\alpha)<\infty$}
Given   $ \theta \in\mathbb{R}$ and  $\epsilon_0>0$,  we say $k$ is an $\epsilon_0$-resonance for $\theta$ if
$ \| 2\theta-k\alpha\|_{\mathbb{R}/\mathbb{Z}}\leq e^{-\epsilon_0|k|}$ and
$\| 2\theta-k\alpha\|_{\mathbb{R}/\mathbb{Z}}=\min_{|j|\leq|k|} \| 2\theta-j\alpha\|_{\mathbb{R}/\mathbb{Z}}$.
\par
Clearly, $0\in \mathbb{Z}$ is  an $\epsilon_0$-resonance.
We order the $\epsilon_0$-resonances $0=|n_0|<|n_1|\leq|n_2|\cdots$. We say  $\theta$ is $\epsilon_0$-resonant if the
set of $\epsilon_0$-resonances is infinite. If $\theta$ is non-resonant, with the set of resonances $\{n_0,n_1,\cdots, n_{j_\theta}\}$,
we set $n_{j_\theta+1}=\infty$.
Notice that  if $\| 2\theta-k\alpha\|_{\mathbb{R}/\mathbb{Z}}=0$ for some $k\in \mathbb{Z}$, then $k$ is an resonance for $\theta$, and
 $\theta$ is not  $\epsilon_0$-resonant.
\par
 Below,  $C$ is a large absolute  constant and $c$ is a small
 absolute  constant, which may change through the arguments, even when appear in the same formula. However, their  dependence on other parameters
 will be explicitly indicated. For instance, we denote by  $C(\alpha)$   a large constant depending on $\alpha$.
\par
Before starting our main work in this part, we firstly give some simple facts.
\begin{lemma} \label{XLe31}
Assume $0<\beta(\alpha)<\infty$, then
\begin{equation}\label{X31}
     \inf_{0<|j|\leq k} ||j\alpha||_{\mathbb{R}/\mathbb{Z}}\geq c(\alpha)e^{-2\beta k},
\end{equation}
and
\begin{equation}\label{X32}
    \inf_{0<|j|\leq k} ||j\alpha||_{\mathbb{R}/\mathbb{Z}}\geq e^{-3\beta k}, \text { for } k>k(\alpha)\footnote{$ k>k(\ast) $ means $k$ is large enough  depending on $\ast$ .}.
\end{equation}
\end{lemma}
\textbf{Proof:}    By $(\ref{G11})$ and  (\ref{XG29}) there exists some $n_0>0$  such that for  $n>n_0(\alpha)$,
\begin{equation}\label{X33}
    ||q_n\alpha||_{\mathbb{R}/\mathbb{Z}}\geq \frac{1}{2}q_{n+1}^{-1}\geq  {e^{-2\beta q_n}}.
\end{equation}
Let $c(\alpha)=\inf_{0<|j|\leq q_{n_0+1}}||j\alpha||_{\mathbb{R}/\mathbb{Z}}>0$.
 Assume   $0<|j|\leq k$.  If $|j|\geq q_{n_0+1}$, select  $q_n\leq |j|<q_{n+1}$ with $n\geq n_0+1$. By (\ref{XG28}) and $(\ref{X33})$
 \begin{eqnarray}
 \nonumber
    ||j\alpha||_{\mathbb{R}/\mathbb{Z}} & \geq& ||q_n\alpha||_{\mathbb{R}/\mathbb{Z}}  \\
   \nonumber
     &\geq& {e^{-2\beta q_n}} \geq  {e^{-2\beta k}} \\
     &\geq& c(\alpha) {e^{-2\beta k}} . \label{X34}
 \end{eqnarray}
If $|j|<q_{n_0+1}$, by the definition of $c(\alpha)$, $$ ||j\alpha||_{\mathbb{R}/\mathbb{Z}}  \geq c(\alpha)\geq c(\alpha) {e^{-2\beta k}}. $$
This implies $(\ref{X31})$.
For $(\ref{X32})$, notice that  $ c(\alpha)>e^{-\beta k}$ for $k>k(\alpha)$. $\qed$
\begin{remark}
In particular, $  ||k\alpha||_{\mathbb{R}/\mathbb{Z}}\geq c(\alpha)e^{-2\beta | k|} $ for  all $k\in \mathbb{Z}\backslash\{ 0\}$. This is a small divisor condition when we solve the
 homological equation (see Theorem \ref{Th417} or Theorem \ref{Th51}).
\end{remark}
\begin{lemma}\label{XXXLe33}
If $\epsilon_0=C_1 \beta>0$, $C_1$ is a  large absolute constant. Then there exists $k_0(\alpha)>0$ such that if $|k|>k_0(\alpha)$  and $||2\theta-k\alpha||\leq e^{-\epsilon_0 |k|}$, then $k$ is  an $\epsilon_0$-resonance for $\theta$.
\end{lemma}
\textbf{Proof:} It suffices to  prove $\| 2\theta-k\alpha\|_{\mathbb{R}/\mathbb{Z}}=\min_{|j|\leq|k|} \| 2\theta-j\alpha\|_{\mathbb{R}/\mathbb{Z}}$.
If $|j|\leq |k|$ and $j\neq k$, by $ (\ref{X32})$ there exists some $k_0(\alpha)$ such that
\begin{eqnarray}
\nonumber
  \| 2\theta-j\alpha\|_{\mathbb{R}/\mathbb{Z}} &\geq&  ||(k-j)\alpha||_{\mathbb{R}/\mathbb{Z}}-||2\theta-k\alpha||_{\mathbb{R}/\mathbb{Z}} \\
   \nonumber
    &\geq&     e^{- 6\beta |k|}-e^{-\epsilon_0 |k|}\\
 &>& e^{-\epsilon_0| k|}\geq \| 2\theta-k\alpha\|_{\mathbb{R}/\mathbb{Z}}\label{XG318}
\end{eqnarray}
for $ k>k_0(\alpha)$. It follows that   $k$ is  an $\epsilon_0$-resonance for $\theta$. $\qed$
\begin{definition}\label{Def31}
 We say that  $ \hat{H}_{\lambda  ,\alpha,\theta} $  satisfies a strong localization
estimate if there exists $ C_0>0$, $ \epsilon_0>0$ and $\epsilon_1>0$  such that for  any solution
 $\hat{H}_{\lambda  ,\alpha,\theta}\hat{u}=E\hat{u}$ with  $\hat{u}_0=1$ and $|\hat{u}_k|\leq 1+|k|$, where $E$ in the spectrum
 of $ \hat{H}_{\lambda  ,\alpha,\theta}$, i.e.,
 $E\in \Sigma_{\lambda,\alpha}$,    we have $|\hat{u}_k|\leq C(\hat{u})e^{-\epsilon_1|k|}$  for $C_0|n_j|<|k|<C_0^{-1}|n_{j+1}|$.
\end{definition}
\par
 \begin{lemma}$ (\text{Lemma } 9.7, \cite{AJ1})$\label{Le31}
Let $\alpha\in \mathbb{R}\backslash \mathbb{Q}$, $x\in\mathbb{R}$ and $0\leq \ell_0 \leq q_n-1$ be such that
$ |\sin\pi(x+\ell_0\alpha)|=\inf_{0\leq\ell\leq q_n-1}    |\sin\pi(x+\ell \alpha)|$, then for some absolute constant $C > 0$,
\begin{equation}\label{G31}
    -Cq_n\leq \sum _{\ell=0,\ell\neq \ell_0}^{q_n-1} \ln|\sin\pi(x+\ell\alpha )|+(q_n-1)\ln2\leq  Cq_n,
\end{equation}
where $  {q_n} $ is given in $ \S 2.5$.
\end{lemma}
The next theorem is our main work in this section.
\begin{theorem}\label{Th32}
   Fix $  \epsilon_0  = C_1\beta >0 $,   where $C_1$   is large enough so that it is   much larger than any absolute constant   $C$, $c^{-1} $ emerging in the present paper.
   Then there exists some constant $C_2$ such that,  for
   $ 0<|\lambda| <e^{-C_2\beta}$,
    $\hat{H}_{\lambda  ,\alpha,\theta}$
  satisfies a strong localization estimate with parameters $C_0=3$, $\epsilon_0=C_1\beta$  and $\epsilon_1= \frac{-\ln|\lambda|}{64}$.
\end{theorem}
\begin{remark}
Refering to  Lemma \ref{XLe42} in the next section, it follows that $|n_{j+1}|>\frac{C_1}{8}|n_j|$. Thus
there exists $k$ such that $3|n_j|<|k|<\frac{1} {3}|n_{j+1}|$ if $C_1$ is large enough.
\end{remark}
 By Aubry  duality
  $\hat{H}_{\lambda  ,\alpha,\theta}=\lambda  {H}_{\lambda^{-1}  ,\alpha,\theta}$, thus  to prove Theorem $ \ref{Th32}$,
  we  only  need  prove $\check{H}_{\lambda   ,\alpha,\theta}\triangleq{H}_{\lambda^{-1}  ,\alpha,\theta} $ satisfies the strong localization estimate
  instead.
 Since this does not   change  any of the statements,  sometimes the dependence of parameters $E,\lambda,\alpha, \theta$ will be ignored in the following.
   Assume $ \check{H}\phi=E\phi$ with $\phi(0)=1$ and $|\phi(k)|\leq1+|k|$. Our objective is to   prove  $|\phi(y)|\leq C(\phi) e^{-\frac{L}{64}|y|}$. Without loss of generality, assume $0<\lambda<1$
    ( for $\lambda<0$, notice that   $\check{H}_{\lambda,\alpha,\theta}=\check{H}_{-\lambda,\alpha,\theta+ \frac{1}{2} }$). By Theorem
   $\ref{Th21}$,  the Lyapunov exponent of  $S_{\lambda^{-1},E}$ satisfies $L=-\ln\lambda$, where $E\in \Sigma_{\lambda^{-1},\alpha}$.
\par
Define $H_I=R_I\check{H}R_I$, where $R_I=$ coordinate restriction to $I=[x_1,x_2]\subset\mathbb{Z}$, and denote by
$ {G}_I =(\check{H}_I-E)^{-1}$  the associated Green  function, if   $\check{H}_I-E$  is invertible. Denote by $ {G}_{I}(x,y)$    the matrix elements of   Green  function
${G}_I$.
\begin{definition}
Fix $m > 0$ and $ 1/10 <\delta<1/2$. A point $y\in\mathbb{Z}$ will be called $(m,k)$-regular with $\delta$ if there exists an
interval $[x_1,x_2]$  containing $y$,  where $x_2=x_1+k-1$ such that
\begin{equation}\label{G32}
  | G_{[x_1,x_2]}(y,x_i)|<e^{-m|y-x_i|} \text{ and dist}(y,x_i)\geq \delta k \text{ for }i=1,2;
\end{equation}
otherwise, $y$ will be called $(m,k)$-singular with $\delta$.
\end{definition}
\par
It is  easy to check that (p. 61, $\cite{B2}$)
 \begin{equation}\label{G33}
   \phi(x)= -G_{[x_1 ,x_2]}(x_1,x ) \phi(x_1-1)-G_{[x_1 ,x_2]}(x,x_2) \phi(x_2+1),
 \end{equation}
 where  $ x\in I=[x_1,x_2] \subset \mathbb{Z}$.
 \begin{lemma}\label{Regular}
For any $m>0$ and any $\delta$ with  $ 1/10 <\delta<1/2$, 0 is $ (m,k)$-singular with $\delta$ if $k>k(m) $.
 \end{lemma}
\textbf{ Proof:} Otherwise,  0 is $ (m,k)$-regular with  some $ 1/10 <\delta<1/2$,
i.e.,\begin{equation}\label{XXG32}
  | G_{[x_1,x_2]}(0,x_i)|<e^{-m|y-x_i|} \leq e^{-\frac{m}{10}k} \text{ for }i=1,2,
\end{equation}
since $|y-x_i|>\frac{k}{10}$.
 In (\ref{G33}), let $x=0$ and  recall that $\phi(x_1-1)\leq 1+|x_1-1|\leq 1+k$, $\phi(x_2+1)\leq 1+|x_2+1|\leq 1+k$.
 Thus
  \begin{equation}\label{XXG33}
   |\phi(0)|= |G_{[x_1 ,x_2]}(x_1,0 ) \phi(x_1-1)+G_{[x_1 ,x_2]}(0,x_2) \phi(x_2+1)|\leq 2 (1+k)e^{-\frac{m}{10}k}.
 \end{equation}
 This implies $|\phi(0)|<1$ if $k>k(m) $, which is contradicted to the hypothesis $\phi(0)=1$. $\qed$
\par
Let us denote
$$ P_k(\theta)=\det(R_{[0,k-1]}(\check{H}_{\lambda,\alpha,\theta}-E) R_{[0,k-1]}),$$
and $A=S_{\lambda^{-1},E}$,
then the $k$-step transfer-matrix $A_k(\theta)$  given  by $(\ref{G22})$   can be written as $(p.14, \cite{B2})$

\begin{equation}\label{G34}
 A_{k}(\theta)=
\left(
  \begin{array}{cc}
   P_k(\theta) &- P_{k-1}(\theta+\alpha)\\
    P_{k-1}(\theta) & - P_{k-2}(\theta+\alpha) \\
  \end{array}
\right).
\end{equation}
    By Cramer's rule (p. 15, $\cite{B2}$) for given  $x_1$ and $x_2=x_1+k-1$, with
     $ y\in I=[x_1,x_2] \subset \mathbb{Z}$,  one has
     \begin{eqnarray}
       |G_I(x_1,y)| &=&  \left| \frac{P_{x_2-y}(\theta+(y+1)\alpha)}{P_{k}(\theta+x_1\alpha)}\right|,\label{G35}\\
       |G_I(y,x_2)| &=&\left|\frac{P_{y-x_1}(\theta+x_1\alpha)}{P_{k}(\theta+x_1\alpha)} \right|.\label{G36}
     \end{eqnarray}
     The numerators in  ($\ref{G35}$) and ($\ref{G36}$) can be bounded uniformly with  respect  to  $\theta$
     by ($\ref{G24} $) and ($ \ref{G34}$), i.e.,
   for  any $ \varepsilon>0$,
     \begin{equation}\label{G37}
        |P_k(\theta)|\leq||A_k(\theta)||\leq e^{(L +\varepsilon)k} \text{ for  sufficiently    large } k  \text{ and   all } \theta.
    \end{equation}
   In fact,  $ (\ref{G37})$ can  be also  uniform  with respect  to $E\in \Sigma_{\lambda^{-1},\alpha}$  by the compactness of $\Sigma_{\lambda^{-1},\alpha} $ and
    subadditivity of $\ln ||A_k||$ (see the proof of  Theorem $\ref{Th46}$ ).
    \par
     Following  $ {\cite{JKS}}$,
     $P_k(\theta)$ is an even function of $ \theta+\frac{1}{2}(k-1)\alpha$  and can be written as a polynomial
      of order $k$ in $\cos2\pi (\theta+\frac{1}{2}(k-1)\alpha )$:

     \begin{equation}\label{G38}
      P_k(\theta)=\sum _{j=0}^{k}c_j\cos^j2\pi (\theta+\frac{1}{2}(k-1)\alpha)    \triangleq  Q_k(\cos2\pi  (\theta+\frac{1}{2}(k-1)\alpha)).
     \end{equation}

 Let $A_{k,r}=\{\theta\in\mathbb{R} \;|\;Q_k(\cos2\pi   \theta   )|\leq e^{(k+1)r}\} $ with $k\in \mathbb{N}$ and $r>0$.
       \begin{definition}
     We  say that the set $\{\theta_1, \cdots ,\theta_{k+1}\}$ is $ \gamma$-uniform if
      \begin{equation}\label{G39}
        \max_{ x\in[-1,1]}\max_{i=1,\cdots,k+1}\prod_{ j=1 , j\neq i }^{k+1}\frac{|x-\cos2\pi\theta_j|}
        {|\cos2\pi\theta_i-\cos2\pi\theta_j|}<e^{k\gamma}.
      \end{equation}
       The next   two lemmas    are from   $\cite{AJ1}$, for self-contain we give the proof.
     \end{definition}
      \begin{lemma}\label{Le33}$(\text{Lemma 9.2 },\cite{AJ1})$
      Suppose $y\in\mathbb{Z}$ is $(L-\rho,k)$-singular with  $1/10<\delta<1/2$, then for any   $  \varepsilon >0$   and   any
      $x\in \mathbb{Z}$ such that $ y-(1-\delta)k \leq x \leq y- \delta k$, we have that $  \theta+(x+\frac{1}{2}(k-1))\alpha   $
      belongs to $A_{k,L-\rho\delta+\varepsilon }$ for    $k$  large enough.
      \end{lemma}
\textbf{ Proof:} Otherwise, there exist   $\varepsilon>0$ and $x_1$ satisfying $ y-(1-\delta)k \leq x_1 \leq y- \delta k$ and
 $   \theta+(x_1+\frac{1}{2}(k-1))\alpha    \notin A_{k,L-\rho\delta+\varepsilon }$, i.e., $P_k(\theta+x_1\alpha)> e^{(k+1)(L-\rho\delta+\varepsilon )}$
 by $(\ref{G38})$.  Let $I=[x_1,x_2]$ with $x_2=x_1+k-1$, then $ y\in I \text{ and }dist(y,x_i)\geq \delta k \text{ for }i=1,2$.
  By $(\ref{G35})$, $(\ref{G36})$ and $(\ref{G37})$, we have
 \begin{equation}\label{1}
    |G_I (y,x_i)|\leq e^{(L+\varepsilon)(k-|y-x_i| )-(k+1)(L-\rho \delta+\varepsilon)}< e^{-(L-\rho)|y-x_i|}\text{ for } i=1,2.
 \end{equation}
 This implies $y$ is $(L-\rho,k)$-regular, contradicting to the hypothesis. $\qed$
      \begin{lemma}\label{Le34}$(\text{Lemma 9.3 },\cite{AJ1})$
      Let $ \gamma_1<\gamma$. If $\theta_1,\cdots,\theta_{k+1}\in A_{k,L-\gamma}$,  then $\{\theta_1, \cdots ,\theta_{k+1}\}$ is not $ \gamma_1$-uniform  for $ k>k(\gamma,\gamma_1,\lambda) $.
      \end{lemma}
 \textbf{ Proof:}
     Otherwise,  $\max_{ x\in[-1,1]} \prod_{ j=1 , j\neq i }^{k+1}\frac{|x-\cos2\pi\theta_j|}
        {|\cos2\pi\theta_i-\cos2\pi\theta_j|}<e^{k\gamma_1} $, $i=1,2,\cdots,k+1$. By $(\ref{G38})$,  we can write polynomial $Q_k(x)$ in the Lagrange
     interpolation form at points $ \cos2\pi  \theta_i$, $i=1,2,\cdots,k+1$, thus
     \begin{eqnarray*}
         |Q_k(x)| &=& \left| \sum_{i=1}^{k+1} Q_k(\cos2\pi\theta_i)\frac{\prod_{j\neq i }(x-\cos2\pi\theta_{j})}{\prod_{j\neq i}(\cos2\pi\theta_i-\cos2\pi\theta_{j})}\right| \\
        &\leq&  (k+1)e^{(k+1)(L -\gamma)} e^{k\gamma_1}=e^{kL}(k+1)e^{-k(\gamma-\gamma_1)+L-\gamma}\\
         &<& e^{kL}
     \end{eqnarray*}
     for all $x\in[-1,1]$     and $ k>k(\gamma,\gamma_1,\lambda) $. By $(\ref{G38})$ again, $|P_k(x)|< e^{kL}$ for  all $ x\in \mathbb{R}$.
    However,  by Herman's subharmonic function methods (see  p.16 $  \cite{B2}$, or p.461 $\cite{Her} $ ), $ \int_{\mathbb{R}/ \mathbb{Z}} \ln|P_{k}(x)|dx\geq k L$.  This is impossible. $\qed$

       Without loss of generality, assume   $3|n_j|<y<\frac{|n_{j+1}|}{3}$. Select $n$ such that  $q_n\leq \frac{y}{8}<q_{n+1}$ and
       let
       $s $ be the largest positive integer satisfying
$ sq_n\leq\frac{y}{8}$. Set $I_1,I_2\subset\mathbb{Z}$ as follows
 \begin{equation}\label{G310}
    I_1=[-2sq_n+1,0] \;and \;I_2=[y-2sq_n+1,y+2sq_n] , \text{ if }n_j<0,
 \end{equation}
\begin{equation}\label{G311}
      I_1=[0, 2sq_n-1]\; and \;I_2=[y-2sq_n+1,y+2sq_n] ,  \text{ if } n_j\geq0.
\end{equation}
In either case, the total number of elements in $I_1  \cup I_2$ is $6sq_n$. Let  $\theta_ {j'} = \theta + j'\alpha$ for $ j' \in I_1  \cup I_2$.
\begin{lemma}\label{Le35}
 Under the  condition   of Theorem $\ref{Th32}$, the set $\{\theta_ {j'}\}_{   j' \in I_1  \cup I_2}$ constructed as $(\ref{G310})$ or $(\ref{G311})$ is
 $C\epsilon_0$-uniform for $y>y(\alpha )$ (or equivalently $n>n(\alpha)$).
\end{lemma}
\textbf{Proof:}
 Firstly we  estimate the numerator in ($\ref{G39}$).  In ($\ref{G39}$), let $x=\cos2\pi a$ and take the logarithm. We have
   $$ \sum _{j' \in I_1  \cup I_2,j'\neq i}\ln|\cos2\pi a-\cos2\pi \theta_{j'}| \;\;\;\;\;\;\;\;\;\;\;\;\;\;\;\;\;\;\;\;\;\;\;\;\;\;\;\;\;\;\;\;\;\;\;\;\;\;\;\;\;\;\;\;\;\;\;\;\;\;\;\;\;\;\;\;$$
$$\;\;\;\;\;\;\;\;\;\;\;\;\;\;\;\;\;\;\;\;\;=\sum_{j' \in I_1  \cup I_2,j'\neq i}\ln|\sin\pi(a+\theta_{j'})|+\sum_{j '\in I_1  \cup I_2,j'\neq i}\ln |\sin\pi(a-\theta_{j'})|
+(6sq_n-1)\ln2  $$
\begin{equation}\label{XG318}
    =\Sigma_{+}+\Sigma_-+(6sq_n-1)\ln2,  \;\;\;\;\;\;\;\;\;\;\;\;\;\;\;\;\;\;\;\;\; \;\;\;\;\;\;\;\;\;\;\;\;\;\;\;\;\;\;\;\;\;\;\;\; \;
\end{equation}
 where
 \begin{equation}\label{XG320}
   \Sigma_{+}=\sum_{j' \in I_1  \cup I_2,j'\neq i}\ln |\sin\pi(a+  \theta_{j'}  )|,
 \end{equation}
 and
  \begin{equation}\label{XG319}
     \Sigma_-=\sum_{j' \in I_1  \cup I_2,j'\neq i}\ln |\sin\pi ( a-\theta_{j'})|.
 \end{equation}
Both $\Sigma_+$ and $\Sigma_-$ consist of $6s$ terms of the form of $(\ref{G31})$, plus 6s terms of the form
\begin{equation}\label{XG321}
    \ln\min_{j'=0,1,\cdots,q_n-1}|\sin\pi(x+j'\alpha)|,
\end{equation}
minus $\ln|\sin\pi(a\pm\theta_i)|$. Since     there exists a interval of
length $q_n$ in sum of (\ref{XG320}) and  (\ref{XG319}) containing $i$, thus the minimum  over this interval is not
 more than $\ln|\sin\pi(a\pm\theta_i)|$ (by the minimality).
       Thus,  by $(\ref{G31})$ one has
 \begin{equation}\label{G312}
   \sum_{j' \in I_1  \cup I_2,j'\neq i}\ln|\cos2\pi a-\cos2\pi \theta_{j'}|\leq-6sq_n\ln2+Cs\ln q_n.
\end{equation}
\par
 The estimate of  the denominator of  ($\ref{G39}$) requires a bit more work. In $(\ref{XG318})$, let $a=\theta_i$,
we obtain
\begin{equation}\label{G313}
    \sum _{j' \in I_1  \cup I_2,j'\neq i}\ln|\cos2\pi \theta_i-\cos2\pi \theta_{j'}| =\Sigma_{+}+\Sigma_-+(6sq_n-1)\ln2,
\end{equation}
 where
 \begin{equation}\label{G315}
   \Sigma_{+}=\sum_{j' \in I_1  \cup I_2,j'\neq i}\ln |\sin\pi(2\theta+ (i+j') \alpha)|,
 \end{equation}
 and
\begin{equation}\label{G314}
     \Sigma_-=\sum_{j' \in I_1  \cup I_2,j'\neq i}\ln |\sin\pi( i-j')\alpha|.
 \end{equation}
Firstly, $\Sigma_-$ consists of $6s$ terms of the form   of  $ (\ref{G31}) $ plus $6s-1$ minimum terms like (\ref{XG321}) (since there exists a interval of
length $q_n$ containing $i$, the sum over this interval is exactly of  the form   $ (\ref{G31}) $).
By  $(\ref{XG28})$ and $(\ref{X33})$,
\begin{equation*}
    \min_{0<|j'| < q_{n+1} }||j'\alpha||_{\mathbb{R}/\mathbb{Z}}  = ||q_n\alpha||_{\mathbb{R}/\mathbb{Z}}  \geq   e^{-2\beta q_n},
\end{equation*}
 for $n>n(\alpha)$.
 Therefore, for $n>n(\alpha)$,
\begin{equation}\label{G316}
 \max\{\ln|\sin x|, \ln|\sin (x+\pi j'\alpha)|\}\geq -C\beta q_n, \text{ for } x\in \mathbb{R} \text{ and } 0 <|j'|<q_{n+1}.
\end{equation}
By known condition $sq_n<q_{n+1}$,  then there exist  at most 6 minimum terms  smaller than $-C\beta q_n $.
Next we estimate the minimum terms. Obviously,   $|i-j'|<Csq_n$ for $i,j' \in I_1\cup I_2$. By $(\ref{X32})$,
\begin{equation}\label{XG327}
     \min _{j' \in I_1  \cup I_2,j'\neq i}\ln|\sin\pi(i-j')\alpha|\geq -Csq_n\beta  \text { for }  n >n(\alpha).
\end{equation}
By   (\ref{G31}),       (\ref{G316}) and (\ref{XG327}),   we obtain
\begin{equation}\label{G317}
     \Sigma_-\geq-6sq_n\ln2-Csq_n\beta.
\end{equation}
\par
 Similarly,  $ \Sigma_{+}$ consist of $6s$ terms of the form of $(\ref{G31})$  plus 6s  minimum terms and
minus $\ln|\sin2\pi \theta_i |$, and  there exist  at most 6 minimum terms   smaller than $-C\beta q_n $ by  $ (\ref{G316})$.
Thus we only need  estimate the minimum term. By  the definition of $I_1$ and $I_2$,
one easily verifies   $i+j'\neq -n_j$ and $|i + j'| < |n_{j+1}|$. By   Lemma \ref{XXXXXLe311} below, one  has
\begin{equation}\label{G319}
      \min _{j' \in I_1  \cup I_2,j'\neq i} ||2\theta+(i+j')\alpha||_{\mathbb{R}/\mathbb{Z}} \geq   e^{-Csq_n\epsilon_0}.
\end{equation}
Replacing (\ref{XG327}) with (\ref{G319}), and following the discussion of $\Sigma_-$, we have
\begin{equation}\label{G320}
     \Sigma_+\geq-6sq_n\ln2-Csq_n\epsilon_0,
\end{equation}
for  $n>n(\alpha )$ or  $y>y(\alpha )$.
Putting $(\ref{G312})$, $(\ref{G317})$ and $(\ref{G320})$   together,
 \begin{equation}\label{G321}
        \max_{ i\in I_1\cup I_2} \prod_{j'\in I_1\cup I_2,j'\neq i } \frac{|x-\cos2\pi\theta_{j'}|}
        {|\cos2\pi\theta_i-\cos2\pi\theta_{j'}|}<e^{C6sq_n\epsilon_0 },
      \end{equation}
      for $y>y(\alpha)$. $\qed$
\begin{lemma}\label{XXXXXLe311}
Under the condition of Lemma \ref{Le35}, suppose  $i+j'\neq -n_j$ and $|i + j'| < |n_{j+1}|$, where   $ i,j' \in I_1\cup I_2$,
then
\begin{equation}\label{XXG319}
      ||2\theta+(i+j')\alpha||_{\mathbb{R}/\mathbb{Z}} \geq   e^{-Csq_n\epsilon_0},
\end{equation}
for $n>n(\alpha)$ (or equivalently   $y>y(\alpha)$ ).
\end{lemma}
\textbf{Proof:} Let   $|k_0|\leq |i+j'|$  be such that $ ||2\theta+k_0\alpha||_{\mathbb{R}/\mathbb{Z}} =\min_{|k|\leq |i+j'|} ||2\theta+k\alpha||_{\mathbb{R}/\mathbb{Z}} $.
\par
Case 1:  $k_0\neq i+j'$.
 If $ ||2\theta+k_0\alpha||_{\mathbb{R}/\mathbb{Z}} \geq   e^{-Csq_n\epsilon_0}$, by the  minimality of $k_0$,
 we have
 \begin{equation*}
    ||2\theta+(i+j')\alpha||_{\mathbb{R}/\mathbb{Z}}\geq ||2\theta+k_0\alpha||_{\mathbb{R}/\mathbb{Z}} \geq   e^{-Csq_n\epsilon_0}.
\end{equation*}
If $ ||2\theta+k_0\alpha||_{\mathbb{R}/\mathbb{Z}} \leq   e^{-Csq_n\epsilon_0}$,    by (\ref{X31})
\begin{eqnarray}
\nonumber
  ||2\theta+(i+j')\alpha||_{\mathbb{R}/\mathbb{Z}} &\geq&  ||(i+j'-k_0)\alpha||_{\mathbb{R}/\mathbb{Z}}-||2\theta+k_0 \alpha||_{\mathbb{R}/\mathbb{Z}}\\
   \nonumber
    &\geq&    c(\alpha) e^{- 2\beta |i+j'-k_0|}-  e^{-Csq_n\epsilon_0}\\
 &\geq&  e^{- Csq_n\beta},\label{XXG318}
\end{eqnarray}
for $n>n(\alpha)$, since $ |i+j'-k_0|<Csq_n$.
\par
Case 2: $k_0=i+j'$. If  $-k_0$ is not  an resonance for $\theta$,  then by the  definition of resonance
\begin{equation*}
    ||2\theta+(i+j')\alpha||_{\mathbb{R}/\mathbb{Z}}\geq  e^{-\epsilon_0|k_0|} \geq   e^{-Csq_n\epsilon_0}.
\end{equation*}
If    $-k_0$ is   an resonance for $\theta$, therefore $|k_0|\leq|n_j|$ (otherwise $-k_0=n_{j+1}$).
Next we discuss $ ||2\theta-n_j\alpha||_{\mathbb{R}/\mathbb{Z}} \geq   e^{-Csq_n\epsilon_0}$ and
$ ||2\theta-n_j\alpha||_{\mathbb{R}/\mathbb{Z}} \leq   e^{-Csq_n\epsilon_0}$ respectively. Following the proof
of case 1,  we also have, for $n>n(\alpha)$
\begin{equation}\label{XX}
      ||2\theta+(i+j')\alpha||_{\mathbb{R}/\mathbb{Z}} \geq   e^{-Csq_n\epsilon_0}.
\end{equation}
 Putting all cases together, we complete the proof of  this lemma.
\par
\begin{remark}
Note that (\ref{XXG319}) holds if $n$ is large enough, which only  depends on $\alpha$, does not depend on  $\theta$.
By the way, all estimates in the present paper  is uniform with respect to $\theta$ and $E\in \Sigma_{\lambda,\alpha}$. This is important.
\end{remark}
 By Lemma  $\ref{Le34} $ and $ \ref{Le35}$, there exists  at least one of $\theta_{j_0}$ with $j_0\in I_1\cup I_2$  such that
   $\theta_{j_0}\notin A_{6sq_n-1,L-C\epsilon_0}$. We will prove that for  all $j'\in I_1$,  $\theta_{j'} \in A_{6sq_n-1,L-C\epsilon_0}$ if
$\lambda<e^{-C_2\beta}$ with $C_2$ large enough, thus there exists some  $j_0\in I_2$  such that $\theta_{j_0}\notin A_{6sq_n-1,L-C\epsilon_0}$.
\begin{lemma}
There exists some  absolute constant $C_2$ such that for all $j'\in I_1$,  $\theta_{j'} \in A_{6sq_n-1,L-C\epsilon_0}$ if $0<\lambda<e^{-C_2\beta}$ and
$n>n(\lambda,\alpha)$.
\end{lemma}
\textbf{Proof:}  Recall that  by Lemma \ref{Regular}, $y=0$ is  $(m, k)$-singular  with any $\delta$ satisfying $\frac{1}{10}<\delta<\frac{1}{2}$ if $k$ is large enough.
 In Lemma $ \ref{Le33}$, let $y=0$, $\delta=\frac{99}{600}$, $\rho= \frac{99}{100}L$, $\varepsilon =\frac{1}{100}L$ and $k=6sq_{n-1}-1$. One easily checks  that
 for all $j'\in I_1$, $\theta_{j'}\in A_{6sq_n-1, \frac{50799}{60000}L}$.  Obviously, $ \frac{50799}{60000}L<L-C\epsilon_0$ if  $0<\lambda<e^{-C_2\beta}$
 because of $L=-\ln \lambda$.   $\qed$

We can now finish \textbf{the proof of Theorem $\ref{Th32}$}.  Let $j_0\in I_2$ be such that $\theta_{j_0}\notin A_{6sq_n-1,L-C\epsilon_0}$.
  Set $I=[j_0-3sq_n+1,j_0+3sq_n-1]=[x_1,x_2]$. Let $\varepsilon=\epsilon_0$  in $(\ref{G37})$,
  combining with $(\ref{G35})$, $(\ref{G36})$, we have
$$|G_I(y,x_i)|\leq e^{(L+\epsilon_0)(6s q_n-1-|y-x_i|)-6sq_n(L -C\epsilon_0)}\leq e^{-L|y-x_i|+Csq_n\epsilon_0}\text{ for }i=1,2.$$
By a simple computation $|y-x_i|\geq sq_n-2\geq \frac{y}{16} $.  Recall that $L=-\ln\lambda$, thus

 \begin{equation}\label{G322}
 |G_I(y,x_i)|\leq e^{- \frac{y}{16} ( L- C \epsilon_0)}\leq e^{- \frac{L}{32}   y }\text{ for }i=1,2,
 \end{equation}
   if  $|\lambda| <e^{-C_2\beta}$   with
 $C_2$  large enough.
 By ($ \ref{G33}$), we obtain that for $ y >y(\lambda,\alpha)$,   $ |\phi (y) |\leq e^{- \frac{L}{64}   y }$  with $y$ satisfying $  3{|n_j|} <  y <  | n_{j+1}|/3$.  This implies  $|\phi (y) |\leq C(\lambda,\alpha)e^{- \frac{L}{64}y}  $  for  all $y$ with $  3{|n_j|} <  y <  | n_{j+1}|/3$.
 For $y<0$, the proof is similar.
 \par
 We   actually have proved a slightly more precise version of Theorem $\ref{Th32}$.
\begin{theorem}\label{XTh38}
 Let $\epsilon_0=C_1\beta$  and $|\lambda|\in (0,e^{-C_2\beta})$ where $C_1,C_2$ are the constants in Theorem $ \ref{Th32}$, and let $\hat u$ be a solution  of the equation $\hat{H}_{\lambda,\alpha,\theta}\hat{u}=E\hat{u}$ satisfying $\hat{u}_0=1$ and $|\hat{u}_k|\leq1+|k|$, where $E\in \Sigma_{\lambda,\alpha}$.  Then we have that $|\hat{u}_k|\leq e^{-\frac{L}{64}|k|}$   if   $3 |n_j|<|k|<3 ^{-1}|n_{j+1}|$  and    $|k|>C(\lambda,\alpha)$, or equivalently,  that $|\hat{u}_k|\leq C(\lambda, \alpha) e^{-\frac{L}{64}|k|}$  for  all $k$ satisfying $3 |n_j|<|k|<3 ^{-1}|n_{j+1}|$, where set $\{n_j\}$ is the $\epsilon_0$-resonance for $\theta$.
\end{theorem}
\begin{remark}\label{xxRe315}
If $\theta$ is not $\epsilon_0$-resonant, and   a solution  $\hat{H}_{\lambda,\alpha,\theta}\hat{u}=E\hat{u}$ satisfying $\hat{u}_0=1$ and $|\hat{u}_k|\leq1+|k|$,  then by Theorem \ref{XTh38},  $|\hat{u}_k|\leq C(\lambda, \alpha) e^{-\frac{L}{64}|k|}$  with  $   |k|>3 |n_{j_\theta}|$,
since $  n_{j_\theta+1} =\infty$, where $L=-\ln \lambda$.
\end{remark}

\section{   The estimate of  rotation number for resonant phase  }
  It is well known that   for  almost every $E\in \Sigma_{\lambda,\alpha}$, there exists a solution  $\hat{u}$ of the equation $\hat{H}\hat{u}=E\hat{u}$ with $\hat u_0=1$ and $|\hat u_k|\le (1+|k|)^C$.\color{black} See  for the proof of    continuous-time Schr\"{o}dinger operator.
The proof of  discrete     Schr\"{o}dinger operator is similar, see   \cite{LIU2}.
 Generally,  it does not hold for
every $E \in \Sigma_{\lambda,\alpha}$. Such exclusion is inherent to Gelfand-Maurin Theorem.
   Avila and Jitomirskaya in ${\cite{AJ2}}$  conquer this  difficulty  by changing the phase $\theta$. This is a starting point of
   the quantitative version of Aubry duality.
\begin{lemma}$(\text{Theorem }3.3, \cite{AJ2})$\label{Le41}
 If $E\in \Sigma_{\lambda,\alpha} $, then there exists $ \theta\in \mathbb{R}$ and a bounded solution  of
 $\hat{H}_{\lambda , \alpha,\theta}\hat{u}=E\hat{u}$ with
 $\hat{u}_ 0  =1$ and $|\hat{u}_ k |\leq 1$.
\end{lemma}

 Fix $\alpha$ such that $0<\beta(\alpha)<\infty$, and fix $C_1$ in Theorem $\ref{Th32}$.
Without loss of generality, assume $\lambda>0$.  By Theorem $\ref{Th32}$ or  Theorem $\ref{XTh38}$, there exists an absolute constant
 $C_2$ such that,  for
    $0<\lambda<e^{-C_2\beta}$,  $\hat{H}_{\lambda  ,\alpha,\theta}$
  satisfies a strong localization estimate with parameters    $ \epsilon_0$, $\epsilon_1=2\pi h $
 and $C_0  $,
 where $\epsilon_0=C_1\beta$, $h=C_1\epsilon_0$ and $C_0=3$. This is because $2\pi h<\frac{-\ln \lambda}{64}$ in view of   $0<\lambda<e^{-C_2\beta}$ with $C_2$ large enough.
Given $E \in \Sigma_{\lambda,\alpha}$, let   $\theta = \theta(E)$ and  $\hat{u}_k$  be given by Lemma $ \ref{Le41}$. In this section,
assume  $\theta(E) $ is $\epsilon_0$-resonant with the infinite set of $\epsilon_0$-resonances   $\{n_j\}_{j=1}^{\infty}$.  Let  $||\cdot||$  be the Euclidean norms, and denote
 $||f||_\eta=\sup _{|\Im x|<\eta} ||f(x)| |$,  $||f||_0=\sup _{x\in \mathbb{R}} ||f(x) ||$.
 Below,   unless stated otherwise, set  $A=S_{\lambda,E} =\left(
                                                                               \begin{array}{cc}
                                                                                 E-2\lambda \cos2\pi x   & -1 \\
                                                                                 1 & 0 \\
                                                                               \end{array}
                                                                             \right)
  $ with $0<\lambda<e^{-C_2\beta}$ and $E\in\Sigma_{\lambda,\alpha}$. Notice that in $\S 3$,  $A=S_{\lambda^{-1},E}$.
\begin{lemma} \label{XLe42}
For $|n_j|$ large enough (depending on $\alpha$),
\begin{equation}\label{XG41}
     ||2\theta-n_j\alpha||_{\mathbb{R}/\mathbb{Z}}\geq e^{-8\beta|n_{j+1}|},
\end{equation}
in particular, $|n_{j+1}|>\frac{C_1}{8}|n_j|$.
\end{lemma}
\textbf{Proof:}   By $(\ref{X32})$,
 \begin{eqnarray}
\nonumber
  \| 2\theta-n_j\alpha\|_{\mathbb{R}/\mathbb{Z}} &\geq&  ||(n_{j+1}-n_j)\alpha||_{\mathbb{R}/\mathbb{Z}}-||2\theta-n_{j+1}\alpha||_{\mathbb{R}/\mathbb{Z}} \\
   \nonumber
    &\geq&    e^{- 6\beta |n_{j+1}|}-e^{-\epsilon_0 |n_{j+1}|}\\
 &\geq& e^{-8 \beta | n_{j+1}|}\label{XG42}.
\end{eqnarray}
This implies ( \ref{XG41}). Combining with the fact
  $||2\theta-n_j\alpha||_{\mathbb{R}/\mathbb{Z}}\leq e^{-\epsilon_0|n_{j }|}$,  one has $|n_{j+1}|>\frac{C_1}{8}|n_j|$. $\qed$

\par
We will say that a trigonometrical polynomial $p : \mathbb{R}/\mathbb{Z}  \mapsto \mathbb{C} $ has essential degree at most $k$ if its
Fourier coefficients   outside an interval $I$ of length $k$ (for $I=[a,b]$, $k=b-a$) are vanishing.
\begin{lemma}($\text{Theorem } 6.1, \cite{AJ2}$  )\label{Le42}
Let $ 1\leq r\leq\lfloor q_{n+1}/q_n\rfloor$.
If $p$ has essential degree at most $k=rq_n-1$ and $x_0\in \mathbb{R}/\mathbb{Z}$, then
\begin{equation}\label{XXXGG45}
    \| p\|_0\leq C q_{n+1}^{ Cr }\sup_{0\leq j\leq k}|p(x_0+j\alpha)|.
\end{equation}
\end{lemma}
In the present paper,  under condition $     \beta(\alpha)=\limsup_{n\rightarrow\infty}\frac{\ln q_{n+1}}{q_n}$,  (\ref{XXXGG45})
becomes
\begin{equation}\label{G42}
    \| p\|_0\leq Ce^{Cr\ln q_{n+1}}\sup_{0\leq j\leq k}|p(x_0+j\alpha)|\leq  e^{C\beta k}\sup_{0\leq j\leq k}|p(x_0+j\alpha)|,
\end{equation}
for $n>n(\alpha)$ or equivalently $k>k(\alpha)$.
\par
  For any  $n$ with $9|n_j|<n<\frac{1}{9}|n_{j+1}|$  of the form $n = rq_k -1 < q_{k+1}$
  ( by Lemma \ref{XLe42}, there exists such $n$  if $|n_j|$ is large enough depending on $\alpha$),
let   $u^{I_1}(x)=\sum_{k\in I_1}\hat{u}_ k e^{2\pi i kx}$  with
$I_1=[-[\frac{n}{2}],n-[\frac{n}{2} ]]=[x_1,x_2]$.  Recall that $\hat{u}_ k $ is given by  Lemma $ \ref{Le41}$ and satisfies the estimate in Theorem \ref{XTh38}. Define    $U^{I_1}(x)= \left(
           \begin{array}{c }
             e^{2\pi i \theta } u^{I_1}(x)\\
             u^{I_1}(x-\alpha)\\
           \end{array}
         \right)
  $,
   by direct computation
  \begin{equation} \label{G43}
    AU^{I_1}(x)=e^{2\pi i\theta}U^{I_1}(x+\alpha)+ e^{2\pi i \theta}\left(
                                     \begin{array}{c}
                                     g(x) \\
                                      0 \\
                                     \end{array}
                                   \right),
  \end{equation}
  and the Fourier coefficients of $g(x)$ satisfy
   \begin{equation}\label{G44}
   \hat{ g}_k=\chi_{I_1}(k)(E-2\cos2\pi(\theta+k\alpha)) \hat{u}_ k -\lambda\sum_{j\in\{-1,1\}}\chi_{I_1}(k-j)  \hat{u}_{ k-j},
  \end{equation}
  where  $ \chi_I$ is the characteristic function of $I$. Since $  \hat{H}  \hat{u} = E \hat{u}$, we also have
  \begin{equation}\label{G45}
   -\hat{ g}_k=\chi_{\mathbb{Z}\backslash {I_1}}(k)(E-2\cos2\pi(\theta+k\alpha) )\hat{u}_ k -\lambda\sum_{j\in\{-1,1\}}\chi_{\mathbb{Z}\backslash {I_1}}(k-j) \hat{u}_{k-j}.
  \end{equation}
By $(\ref{G44})$ and $(\ref{G45})$, $\hat{g}_k\neq 0$, only at four  points $x_1$, $x_2$, $x_1-1$ and $x_2+1$. By the  strong localization estimate $|\hat{u}_k|\leq  C(\lambda,\alpha)e^{-2\pi h |k|}$,
 it is easy to see $||g||_\frac{h}{3}\leq   C(\lambda,\alpha) e^{-   3 hn}$, in particular $||g||_\frac{h}{3}\leq    e^{-   2 hn}$ for $n>n(\lambda,\alpha)$, since $C(\lambda,\alpha)<e^{h n}$  for $n>n(\lambda,\alpha)$.
\par
\begin{lemma}$(\text{Theorem }10, \cite{A3})$ \label{Le43}
If $ \alpha\in \mathbb{R}\backslash \mathbb{Q}$,  $\lambda\neq 0$, $E\in  \mathbb{R}$ and $ \epsilon\geq 0$, then
$ L(\alpha, A^  {E ,\epsilon} )= \max\{L(\alpha,A^E),  \ln|\lambda| +2\pi \epsilon\}$,
where $$A^ {E,\epsilon}  =\left(
                                                                               \begin{array}{cc}
                                                                                 E-2\lambda \cos2\pi(x+i\epsilon) & -1 \\
                                                                                 1 & 0 \\
                                                                               \end{array}
                                                                             \right),
 $$
 and $A^E= A^  {E ,0}$.
\end{lemma}
\begin{corollary} \label{Co44}
If $ \alpha\in \mathbb{R}\backslash \mathbb{Q}$,  $|\lambda|<1$ and $ \frac{ \ln|\lambda|}{2\pi}\leq\epsilon\leq\frac{-\ln|\lambda|}{2\pi}$, then
$ L(\alpha, A^ {E,\epsilon}  )=0$ for $E\in \Sigma_{\lambda,\alpha}$.
\end{corollary}
\textbf{Proof:} By Theorem $ \ref{Th21} $, if $|\lambda|<1$ and $E\in \Sigma_{\lambda,\alpha}$, then
$L(\alpha,A^E)=0$. If $0\leq\epsilon\leq \frac{-\ln|\lambda|}{2\pi}$, i.e.,
  $ \ln|\lambda| +2\pi \epsilon \leq0 $, then   $ L(\alpha, A^ {E,\epsilon } )=0$ by   Lemma $\ref{Le43}$.  By symmetry $   L(\alpha, A^ {E, \epsilon } )=0$ for $ \frac{ \ln|\lambda|}{2\pi}\leq\epsilon\leq\frac{-\ln|\lambda|}{2\pi}$ . $\qed$
\par
Next we will set up the priori estimate of  transfer matrix,  precisely,  $  ||A_n(x)||=e^{o(n)}$ through band  $|\Im x|< \frac{-\ln|\lambda|}{2\pi}$ uniformly, where $A=S_{\lambda,E}$ and $A_n $ is given by  (\ref{G22}).
This can be done by Furman's uniquely ergodic theorem and  vanishing  Lyapunov exponent (Corollary \ref{Co44}).
\begin{theorem}$(\text{Theorem } 1, \cite{FUR})$\label{Th45}
Let $\{f_n\}$  be a continuous subadditive cocycle on a uniquely
ergodic system $(X, \mu, T)$,
i.e., $X$ is a compact metric space, $ T: X \mapsto X $ is a homeomorphism with
$\mu$ being the unique $T$-invariant probability measure on $X$, and $f_n \in  C(X)$ with $f_{n+m}(x) \leq f_n(x)+f_m(T^nx)$
for all $x\in X$. Then for every $x \in X $ and uniformly on $X$:
\begin{equation}\label{fur}
  \limsup_{n\rightarrow \infty}\frac{1}{n}f_n(x)\leq \lim_{n\rightarrow \infty}\frac{1}{n}\int f_n d\mu.
\end{equation}
\end{theorem}
Theorem \ref{Th45} is  usually called  Furman's uniquely ergodic theorem.
\begin{theorem} \label{Th46}
   $  ||A_k^E(x)||_{\eta}\leq C (\lambda,\alpha)  e^{\beta k}$  for all $E\in \Sigma_{\lambda,\alpha}$,
    where $\eta=\frac{- \ln |\lambda|}{2\pi}$.
\end{theorem}
\textbf{Proof:} By Corollary $ \ref{Co44}$, $ L(\alpha, A^{E, \epsilon } )=0$ for any $-\eta \leq\epsilon \leq \eta$ and $E\in \Sigma_{\lambda,\alpha}$.
In Theorem $ \ref{Th45}$, let $ f_n=\ln ||A_n^{E,\epsilon}|| $, $X=\mathbb{R}/ \mathbb{Z}$, $Tx=x+\alpha$ and $\mu$ is  Lebesgue measure.
Since irrational rotations are uniquely ergodic, then there exists some $k_0(\lambda,\alpha,E,\epsilon)$ such that
\begin{equation*}
    \ln||A_ k^ { E } (x)||<\beta k
\end{equation*}
  for all   $x$   satisfying  $\Im x=\epsilon $  with $|\epsilon|\leq\frac{-\ln|\lambda|}{2\pi}$ and $ k\geq k_0(\lambda,\alpha,E,\epsilon)$.
 By continuity and compactness of $ \mathbb{R}/ \mathbb{Z}   $,  there exists $\delta(\lambda,\alpha,E,\epsilon)$ such that if $ |E'-E|<\delta$ and
$| \Im x' -\epsilon|<\delta$, then
\begin{equation} \label{G47}
    \ln||A_k^{E' }(x')||<\beta k
\end{equation}
for every $k_0(\lambda,\alpha,E,\epsilon)\leq k\leq 2k_0(\lambda,\alpha,E,\epsilon)+1$.
\par
For any  $k>2k_0(\lambda,\alpha,E,\epsilon) $,
let $k=sk_0+r$, where $k_0 \leq r<2k_0$, then by subadditivity,
\begin{equation*}
    \ln||A_k^{E' }(x')||\leq s \max_{| \Im x_1 -\epsilon|<\delta}  \ln ||A_{k_0}^{E' }(x_1)||+\max_{| \Im x_1 -\epsilon|<\delta}  \ln ||A_{r}^{E' }(x_1)|| <\beta k.
\end{equation*}
Thus $(\ref{G47})$ holds for all $k\geq k_0(\lambda,\alpha,E,\epsilon)$.
 By the compactness of $ \{ | \epsilon|\leq \eta\}$ and $\Sigma_{\lambda,\alpha}$, there exists $k_0(\lambda,\alpha)$,
 such that $$ \ln||A_k^{E}(x)||<\beta k $$
   for every $ x $ satisfying  $|\Im x|\leq \eta$, $E\in \Sigma_{\lambda,\alpha}$ and $k>k_0(\lambda,\alpha )$.
   It follows that
    $$  ||A_k^E(x)||_\eta\leq C (\lambda,\alpha)  e^{\beta k}.$$
    We finish the proof.
    \begin{remark}\label{Re47}
   In fact, our proof suggests    that  for any   $\delta>0$, $ ||A_k^E(x)||_\eta\leq C (\delta,  \lambda,\alpha)e^{\delta k}$
    with $ \eta=-\frac{1}{2\pi} \ln|\lambda|$. This verifies a claim by Avila  in the footnote 5 of  $\cite{A2}$.
    \end{remark}
For more subtle estimate of the transfer matrix, a couple of lemmata and theorems are necessary.
 \begin{theorem}$\label{Th48}$
 For $n>n(\lambda,\alpha)$,
  \begin{equation}\label{G48}
    \inf_{|\Im x|<\frac{h}{3}}\|U^{I_1}(x)\|\geq e^{-C \beta n}.
  \end{equation}
  \end{theorem}

 \textbf{ Proof}:
  Otherwise, let $x_0$  with $\Im x_0=t$  and $|t|<\frac{h}{3}$ such that   $ \|U^{I_1}(x_0)\|\leq e^{-C \beta n}$.
  By ($ \ref{G43}$) and  Theorem $\ref{Th46}$, $||U^{I_1}(x_0 +j\alpha)||\leq e^{-C\beta n} $, $0 \leq j\leq n$, since $||g||_{\frac{h}{3}}<e^{-2h n}$ for
   $n>n(\lambda,\alpha)$. This implies for  $n>n(\lambda,\alpha)$, $|u^{I_1}(x_0 +j\alpha)|\leq e^{-C\beta n} $, $0 \leq j\leq n$.
   Thus $ \|u^{I_1}_t\|_0\leq e^{-C\beta  n}$ by  $ (\ref{G42})$, contradicting  to $\int u^{I_1}_t(x)dx=1$ (since $\hat{u}_0=1$), where $u^{I_1}_t(x)=u^{I_1}(x+ti)$. $\qed$
  \begin{theorem}$(\text{Theorem } 2.6,\;  \cite{A2} )$\label{Th49}
  Let $U:\mathbb{R}/ \mathbb{Z}\rightarrow \mathbb{C}^2$ be analytic  in $|\Im x|<\eta$. Assume that
   $ \delta_1<||U(x)|| <\delta^{-1}_2$ holds for all $x$ satisfying $|\Im x|<\eta$ .
  Then there exists  $B:\mathbb{R}/ \mathbb{Z}\rightarrow \text{SL}(2,\mathbb{C} )$  being  analytic  in $|\Im x|<\eta$ with first column $U$ and
  $||B||_\eta\leq C\delta_1^{-2}\delta_2^{-1}(1-\ln(\delta_1\delta_2))$.
  \end{theorem}
  %\begin{remark}
  %\color{red} The Corona Theorem states that if $d \geq 1$ and $a_i: \mathbb{D} \mapsto \mathbb{C}, 1 \leq i \leq d$ are bounded
%holomorphic functions such that $\max_i  |a_i|\geq \delta$ pointwise then there exist bounded
%holomorphic functions $b_i: \mathbb{D} \mapsto \mathbb{C}, 1 \leq i \leq d$ such that  $\sum _{i=1}^da_ib_i=1$, where $\mathbb{D}=\{x\in \mathbb{C}:|x|<1\}$.
 % Trent prove another form of  Theorem \ref{Th49} in $ {\cite{T}}$, i.e.,  the involved   vector and matrix  are  defined in disk $\mathbb{D}$, the version  in annulus
 % $\{x\in \mathbb{C}/ \mathbb{Z}:|\Im x|<\eta\}$ is given in  \cite{A2}.  ????? what's use?\color{black}
  %\end{remark}
  \begin{lemma}\label{Le411}
  \begin{equation}\label{G410}
  \max_{ x \in \mathbb{R}} \|A_m(x)\| \leq  C(\lambda, \alpha) m^C  .
  \end{equation}
  \end{lemma}
  \textbf{Proof:}
  The estimate $|\hat{u}_k |\leq1$ implies
  $||U^{I_1}||_\beta<e^{C\beta n}$.
  Let $B(x)\in \text{SL} (2,\mathbb{C})$ be the matrix,  whose first column is  $U^{I_1}(x)$,   given by  Theorem \ref{Th49}  with
  $\eta=\beta$,    then
$||B||_\beta\leq    e^{C\beta n} $ for $n>n(\lambda,\alpha)$. Combining with  $(\ref{G43})$,  one easily verifies
\begin{equation}\label{G411}
     B(x+\alpha)^{-1}A(x)B(x)=
     \left(
       \begin{array}{cc}
         e^{2\pi i \theta} & 0\\
         0 &  e^{-2\pi i \theta}\\
       \end{array}
     \right)+
     \left(
       \begin{array}{cc}
         \beta_1(x) & b(x) \\
         \beta_2(x) & \beta_3{x} \\
       \end{array}
     \right),
\end{equation}
where $\|b \|_{ \beta}< e^{C\beta n}$, and $\|\beta_1 \|_{ \beta}$, $\|\beta_2 \|_{ \beta}$, $\|\beta_3 \|_{ \beta}<  e^{- h n}$, since
$||g||_\frac{h}{3}\leq    e^{-   2 hn}$.
Taking $ \Phi=DB(x)^{-1}$, where $D=\left(
                                      \begin{array}{cc}
                                        d & 0 \\
                                        0 & d^{-1}\\
                                      \end{array}
                                    \right)
$
 with $d =  e^{- \frac{hn}{4}}$, we get
 \begin{equation}\label{G412}
     \Phi(x+\alpha) A(x)\Phi(x)^{-1}=
     \left(
       \begin{array}{cc}
         e^{2\pi i \theta} & 0\\
         0 &  e^{-2\pi i \theta}\\
       \end{array}
     \right)+
      H(x)
\end{equation}
 where $\|H\|_{ \beta}<  e^{- \frac{hn}{4}}$ and $\|\Phi\|_{ \beta}<e^{ hn}$. Thus
 \begin{equation}\label{G413}
    \sup_{0\leq s \leq e^{\frac{ hn}{4}}}\|A_s(x)\|_{ \beta} \leq  e^{3hn}.
 \end{equation}
 If $m>C(\lambda,\alpha)$, we can select $n$ with  $ C\frac{\ln m}{h}<n<C^2\frac{\ln m}{h}$
   of the form $n = rq_k -1 < q_{k+1}$ and $9|n_j|<n<\frac{1}{9}|n_{j+1}|$,  thus  $||A_m||_\beta<   m^C $  by   $(\ref{G413} )$.
   That is  $||A_m||_\beta<  C(\lambda,\alpha) m^C $ for all $m$. $\qed$
 \par
 Fix some $n = |n_j| $ and let $N = |n_{j+1}| $. Construct  new  function $ u^{I_2}(x)$ with $I_2=[-[\frac{N}{9}], [\frac{N}{9}]]$ and
      vector-valued function  $U^{I_2}(x)= \left(
           \begin{array}{c }
             e^{2\pi i \theta } u^{I_2}(x)\\
             u^{I_2}(x-\alpha)\\
           \end{array}
         \right)
  $ as before.
  \begin{theorem}\label{Th412}
  For $n>n(\lambda,\alpha)$,
  \begin{equation}\label{G414}
    \inf_{|\Im x|<\frac{h}{3}}\|U^{I_2}(x)\|\geq   e^{-C \beta n}.
  \end{equation}
  \end{theorem}
 \textbf{ Proof:} Let  $rq_k$  be the minimal  such that $rq_k>9|n_j|$ and $rq_k-1<q_{k+1}$,
 and let $J = [-[\frac{rq_k}{2} ] ,  rq_k-1 -[\frac{rq_k}{2} ] ]$. Define  $U^J (x)$ as before.
   By the estimates $|\hat{u}_k|\leq   e^{-2\pi h|k|}$ for
 $3n<|k|<\frac{N}{3}$ and $|u_k|\leq 1$ for others (since $ n>n(\lambda,\alpha)$ ),   we have $||U^{I_2} - U^J||_\frac{h}{3}\leq  e^{-  hn  }  $.
  By $ (\ref{G48})$ and a simple fact $rq_k\leq Cn$,
   one has   $\inf_{|\Im x|<  \frac{h}{3}}\|U^J(x)\|\geq   e^{-C\beta n}$.   This implies
\begin{equation}\label{G415}
\inf_{|\Im x|<\frac{ h}{3}}\|U^{I_2}(x)\|\geq  e^{-C\beta n}.
\end{equation}
We finish the proof of the theorem.

\par
Let $\tilde{U}(x)=e^{\pi i n_j x}U^{I_2}(x)$ and $ \tilde{\theta}=\theta-\frac{n_j\alpha}{2}$.
Note that $\tilde{U}(x)$ depends on $I_2$, for simplicity we drop the dependence, since below the interval is always  $I_2=[-[\frac{N}{9}], [\frac{N}{9}]]$.
  Let $B(x)$ be the matrix with columns  $\tilde{U}(x)$ and $\overline{\tilde{U}(x)}$, where  $\overline{\tilde{U}(x)}$ is the complex conjugate of
  $\tilde{U}(x)$,
   and let  $P^{-1}=\|2\theta-n_j\alpha\|_{\mathbb{R}/\mathbb{Z}}$.
   By the same arguments of $(\ref{G43})-(\ref{G45})$,  for $n>n(\lambda,\alpha)$,
   \begin{equation}\label{G416}
    A\tilde{U}(x)=e^{2\pi i \tilde{\theta}}\tilde{U}(x+\alpha)+\left(
                                     \begin{array}{c}
                                      g(x) \\
                                      0 \\
                                     \end{array}
                                   \right)
                                   \text{ with } \|g\|_{ \frac{h}{3} } <e^{- c hN}.
  \end{equation}

  By the  definition of resonance  and Lemma $ \ref{XLe42}$,
  \begin{equation}\label{XG419}
  e^{\epsilon_0 n}\leq P\leq e^{8\beta N } \text{ for   }  n>n( \alpha) .
  \end{equation}

  \begin{theorem}\label{Th414} For $n>n(\lambda,\alpha)$,
  \begin{equation}\label{G417}
    \inf_{x\in \mathbb{R}/\mathbb{Z}}|\det B(x)|\geq    P^{-C}.
  \end{equation}
  \end{theorem}
 \textbf{ Proof}:
 By the proof of Lemma 8.1 in $ \cite{AJ2}$,   for any complex matrix $M$ with columns $V$ and $W$,
 \begin{equation}\label{G418}
    |\det M|= ||V||\min_{\lambda\in\mathbb{ C}}||W-\lambda V||
 \end{equation}
 and the minimizing $\lambda$ satisfies $||\lambda V||\leq ||W||$.
 Suppose   (\ref{G417}) would not hold. By Theorem $\ref{Th412}$  $\inf_{ x\in \mathbb{R}}\|\tilde{U} (x)\|\geq e^{-C\beta n}$,   then there exists
 $x_0\in \mathbb{R}$
 and $\lambda_0\in \mathbb{C}$ ($|\lambda_0|\leq 1$)
 such that  $ |\overline{\tilde{U}(x_0)}-\lambda_0 \tilde{U}(x_0) |\leq    P^{-C}$.
  By $(\ref{G410})$ and $ (\ref{G416})$, we have
  \begin{equation}
      | e^{-2\pi ij\tilde{\theta}} \overline{\tilde{U}(x_0+j\alpha)}- e^{ 2\pi ij\tilde{\theta}} \lambda_0\tilde{U}(x_0+j\alpha)  |
     \leq   P^{-C},0 \leq j \leq  P.
  \end{equation}
  That is
   \begin{equation}\label{G419}
      |   \overline{\tilde{U}(x_0+j\alpha)}- e^{ 4\pi ij\tilde{\theta}} \lambda_0\tilde{U}(x_0+j\alpha)  |
     \leq   P^{-C},0 \leq j \leq  P.
  \end{equation}
Notice a   simple fact $|e^{  4\pi ij\tilde{\theta}}-1|<C||2j \tilde {\theta}||_{\mathbb{R}/\mathbb{Z}}< P^{-c}$, for $0 \leq j \leq  P^ {1-c}$,
since $||2 \tilde{\theta}||_{\mathbb{R}/\mathbb{Z}}= P^{-1}$.
  Combining with (\ref{G419}) and noting $||\tilde{U}||_0\leq C(\lambda,\alpha)n$ by the strong localization estimate,  one has
  \begin{equation}\label{G420}
     |   \overline{\tilde{U}(x_0+j\alpha)}- \lambda_0 \tilde{U}(x_0+j\alpha)  |
     \leq   P^{-c},0 \leq j \leq  P^ {1-c}.
  \end{equation}

 Denote by   $\tilde{U}_k(x)$    truncating the Fourier coefficients  of  $\tilde{U} (x)$  at scale $k=\frac{c}{\beta} \ln P$.
  By $(\ref{XG419})$, one has  $9 n<k< \frac{1}{9}N$. By the strong localization estimate in  Theorem $\ref{XTh38} $ and the definition of $\tilde{U} (x)$,
  \begin{equation}\label{G421}
   \|\tilde{U}-\tilde{U}_k\|_0\leq e^{-\frac{c}{\beta} h\ln P  }\leq P^{-c}.
  \end{equation}
Therefore, we may assume the essential  degree of $\tilde{U}$  is $\frac{c}{\beta}\ln P$.
 By $(\ref{G42})$ and (\ref{G420}),  we have  (first replacing $\tilde{U}(x)$ with $\tilde{U}(2x)$ so that $\tilde{U}(2x)$ is well defined in $\mathbb{R}/\mathbb{Z}$)
\begin{equation}\label{G422}
\sup_{x\in \mathbb{R}/\mathbb{Z}} |\overline{\tilde{U}(x)}- \lambda_0\tilde{U}(x) | \leq e^{C\beta \frac{c}{\beta}\ln P}P^{-c}\leq  P^{-c}.
\end{equation}
In (\ref{G419}), let $j=[\frac{P}{4}]$, we get
\begin{equation}\label{G423}
   |i\overline{\tilde{U}(x_1)}+i\lambda_0 \tilde{U}(x_1) |\leq P^{-c},
\end{equation}
where $x_1=x_0+[\frac{P}{4}]\alpha$. By $(\ref{G422})$ and $(\ref{G423})$,
 $ |\tilde{U}(x_1) |\leq P^{-c}$.  Recall that   $ \inf_{ x \in \mathbb{R}}\|\tilde{U}(x)\|\geq e^{- C\beta n}$, thus
 we get $P \leq e^{C\beta n}$. This  contradicts  to $(\ref{XG419})$   $P\geq e^{\epsilon_0n}$, since  $\epsilon_0=C_1\beta $ and
 we assume
 $C_1$ is much larger than any absolute constant $C$ emerging in this paper. $\qed$
 \par
 The following theorem gives a sharp estimate of the rotation number if phase $ \theta(E)$ is $\epsilon_0$-resonant.
 \begin{theorem}\label{Th415}
  Fix    $n=|n_j| $  (large enough depending on $\lambda$ and $\alpha$ ) and   $N=|n_{j+1}|$,
then there exists  $m_j$ with $|m_j|\leq C n $  such that $ || 2\rho(\alpha,A) -m_j\alpha\pm(2  \theta -n_j \alpha) ||_{\mathbb{R}/\mathbb{Z}}\leq e^{-ch N}$.
  \end{theorem}
 \textbf{ Proof:}
 Let  $S=Re\tilde{U}$, $T=\Im \tilde{U}$, and let  $\tilde{W}$ be the matrix with columns $S$ and $\pm T$
so that  $\det \tilde{W} > 0$. Then by $ (\ref{G416})$
\begin{equation}\label{G424}
    A\tilde{W}(x)= \tilde{W}(x+\alpha)\cdot R_{\mp \tilde{\theta}} +O(e^{- chN}), x\in\mathbb{ R}/\mathbb{Z}.
  \end{equation}
  Let $W(x)= |\frac{\det B(x) }{2} |^{-1/2} \tilde{W}(x)$,  it is easy to verify   $\det W = 1$. By  Theorem $ \ref{Th414}$,
  \begin{equation}\label{G425}
    A W (x)= \frac{|\det B(x+\alpha)|^{1/ 2}}{|\det B(x)|^{1/ 2}}  W(x+\alpha)\cdot R_{\mp \tilde{\theta}} +O(e^{- chN }), x\in\mathbb{ R}/\mathbb{Z}.
  \end{equation}
  By $( \ref{G416})$ and $\det A=1$, $|\det B(x+\alpha) |-|\det B(x) |= O(e^{- chN })$,
  thus we have
  \begin{equation}\label{G426}
     A W (x)=   W(x+\alpha)\cdot R_{\mp \tilde{\theta}} +O(e^{- chN }), x\in\mathbb{ R}/\mathbb{Z}.
  \end{equation}
Since $\det W=1$ and  $W(x)= |\frac{\det B(x) }{2} |^{-1/2}| \tilde{W}(x)$,   $||W^{-1}|| \leq P^{C} \text{ for }  x\in\mathbb{ R}/\mathbb{Z}$.
Then
\begin{equation}\label{G427}
  W(x+\alpha)^{-1}  A W (x)=     R_{\mp \tilde{\theta}} +O(e^{- chN }), x\in\mathbb{ R}/\mathbb{Z}.
\end{equation}
 Since  $  W(x)$   is well defined
in $\mathbb{R}/2\mathbb{Z}$, combing with (\ref{XXG26}),
$$ ||2\rho(\alpha,A)-m \alpha \pm 2\tilde{\theta}| |_{\mathbb{R}/\mathbb{Z}} \leq   e^{ -chN  } , $$
where $ m =\deg(W)$.
 Thus, to prove this theorem, we only need   prove $|\deg(W)|\leq Cn$.
\par
Next we will  estimate the degree of $W$. The  degree of $W$ is the same as the degree of any of its columns\footnote{Let $S: \mathbb{R}/2\mathbb{Z}\mapsto \mathbb{R}^2\backslash 0$, we say degree of $S$ is $k$, denoted by $\deg(S)=k$, if  $S$ is homotopic
to $\left(\begin{array}{c}
               \cos\frac{k}{2}x \\
                \sin \frac{k}{2}x
             \end{array}\right )
$.}.
 It is enough to estimate the degree of
$ \frac{M(x)}{\| M(x)\|}$ for $M=S$or $M=T$. Notice that $ ||\int_{\mathbb{R}/ \mathbb{Z}} e^{-\pi i n_j x}(S(x)+iT(x))dx ||= \sqrt{2} $.
Without loss of generality, assume    $\int_{\mathbb{R}/ \mathbb{Z}} ||S(x) ||dx \geq \sqrt{2}/2$.
By $(\ref{G424})$,
\begin{equation}\label{G428}
   AS(x)=S(x+\alpha)\cos 2\pi\tilde{\theta}  \mp T(x+\alpha)\sin 2\pi\tilde{\theta}  +O(e^{-chN}),x\in\mathbb{ R}/\mathbb{Z}.
\end{equation}
Combining with $||2\tilde{\theta}||_{\mathbb{R}/\mathbb{Z}}\leq e^{-\epsilon_0n}$,
 we have
$A S(x)=   S(x+\alpha)  +O(e^{-  c\epsilon_0n })$, or $A S(x)=  - S(x+\alpha)  +O(e^{-  c\epsilon_0n })$.
Following the proof of Theorem $\ref{Th412}$, we have the similar  estimate
  \begin{equation}\label{G429}
    \inf_{x\in \mathbb{R}}\|S(x)\|\geq   e^{- C\beta n }.
  \end{equation}
 Denote    $\tilde{S}(x)$  by truncating the Fourier series of $S$ at scale $C  n $, then
  $$\|\tilde{S}(x)-S(x)\|\leq    e^{-Ch n}<\frac{\| S(x)\|}{2} $$
  for  $  x\in\mathbb{ R}/2\mathbb{Z}$ and $n>n(\lambda,\alpha)$.
  Thus  the degree of $ {S}$ is equal to the degree of $\tilde{S}$.
  \par
  Now we estimate the degree of $\tilde{S}(x)$. Let $\tilde{S}(2x)=\left (\begin{array}{c}
                                                                               \tilde{S}_1(x) \\
                                                                              \tilde{S}_2(x)
                                                                             \end{array}
                                                                             \right)
  $, then $   \tilde{S}_1(x) $,  $\tilde{S}_2(x)$  only have Fourier series   at scale $Cn$.  Notice that
   $ \tilde{S}_1(x)+i\tilde{S}_2(x)$ can be written as a polynomial of $z$ and $z^{-1}$, where $z=e^{2\pi i x}$.
   More precisely, there exists a polynomial $f(z)$  of order less than  $   Cn$ and $k \in \mathbb{N} $ such that $\frac{f(e^{2\pi i x})}{e^{2\pi i k x}}= \tilde{S}_1(x)+i\tilde{S}_2(x)$,  where  $k<Cn$.
It is a well known fact  that   the  degree  of $\tilde{S}( x)$ is equal to  the zeros of $f(z)$ in disk $\mathbb{D}=\{z:|z|\leq 1\}$ minus $k$. Then
       $|\deg \tilde{S}| \leq Cn $, i.e.,
  $|\deg W| \leq Cn$. $\qed$
  \begin{remark}
  From (\ref{G427}), it is easy to see that $S_{\lambda,E}$ is almost reducible to $R_{\pm \theta}$,  if  $\theta=\theta(E)$
    given by Lemma \ref{Le41} is $\epsilon_0$-resonant. Combining with Theorem \ref{Th417} in the next section, we have
     for every $ E\in \Sigma_{\lambda,\alpha}$, $S_{\lambda,E}$ is almost reducible.
  \end{remark}

\section{ Reducibility for non-resonant phase}
In $\S 4$, we obtain  sharp  estimate of  the rotation number $\rho(\alpha,A)$  when  $\theta(E)$ is   $\epsilon_0$-resonant. In this section, we will
set up   reducibility for $A=S_{\lambda,E}$  with $E\in \Sigma_{\lambda,\alpha}$  when  $\theta(E)$ is not    $\epsilon_0$-resonant.
\begin{lemma} \label {matrix B}
Let $W: \mathbb{R}/2\mathbb{Z}\mapsto  \mathbb{C}^2$ be an real analytic vector   in  $| \Im x|<\eta$.  Assume that $
\inf _ {| \Im x|<\eta}\|W(x)\| >\delta$ with some $\delta>0$,  then there exists $B: \mathbb{R}/2 \mathbb{Z}
\mapsto    \text{SL}(2,\mathbb{ R})$  being real analytic  in $| \Im x|<\eta$  with first column $W$.
\end{lemma}
 \textbf{{Proof:}} Let $W(x)=\left (
 \begin{array}{c}
                            w_1(x) \\
                            w_2(x)
                          \end{array}\right)
$, by Theorem $ \ref{Th49}$ there exist  $b_1$ and $b_2$ being  analytic   in  $| \Im x|<\eta$ such that $w_1b_1-w_2b_2=1$.
Let $\tilde{w}_1(z)=\frac{b_1(z)+\overline{b}_1({z})}{2}$ \footnote{$\overline{a} (z) \text{ is defined by } \overline{a} (z)=\sum \overline{a_n}z^n, \text{ if }
 {a} (z)=\sum  {a_n}z^n $. Notice  that $ \overline{a (z)} $ is the complex conjugate of $a (z)$, however $\overline{a}(z)$ is not. } and $\tilde{w}_2(z)=\frac{b_2(z)+\overline{b}_2( {z})}{2}$,
then  $B=\left(
                                                               \begin{array}{cc}
                                                                 w_1 & \tilde{w}_2\\
                                                                 w_2 & \tilde{w}_1 \\
                                                               \end{array}
                                                             \right):
\mathbb{R}/ 2\mathbb{Z}
\mapsto    \text{SL}(2,\mathbb{ R})$ is  real analytic in $|\Im x|<\eta$.  $\qed$
\begin{remark}\label{XRe59}
 Given  a non-zero real         analytic vector-valued function  $W(x) =\left(\begin{array}{c}
                                        w_1(x) \\
                                        w_2(x)
                                      \end{array} \right)
 $ with $ W(x+1)=\pm W(x)$, all of  Avila, Jitomirskaya, Puig and so on  construct $B$ as follows:  $$B(x)=\frac{1}{w_1^2+w_2^2}\left(
                                                    \begin{array}{cc}
                                                       w_1(x) & -w_2(x) \\
                                                      w_2(x) &  w_1(x)\\
                                                    \end{array}
                                                  \right).
 $$
 Since both $w_1 $ and $ w_2$ are real  analytic,
$ w_1^2+w_2^2>0 $  for $x\in \mathbb{R}$. By continuity, $ w_1^2+w_2^2\neq0 $ in a neighbor of real  axis  and
   $ B :\mathbb{R}/  \mathbb{Z}\mapsto \text{PSL}(2,\mathbb{R})$ is real         analytic in a  neighbor of real axis
 (this process is a key step to
 set up reducibility for cocycle $A=S_{\lambda,E}$. See the proof of Theorem \ref{Th417}). Usually, $ B(x)$ is not real         analytic in
  the given strip.
 In the present paper, since $W(x) $ is well defined in $ \mathbb{R}/  2\mathbb{Z}$, we  can use Lemma \ref{matrix B} to construct a cocycle
 $B: \mathbb{R}/2 \mathbb{Z}
\mapsto    \text{SL}(2,\mathbb{ R})$ with first    column $W$ so that
 $B$ is real         analytic in the given strip. However,
 we    do not have  $B: \mathbb{R}/  \mathbb{Z}
\mapsto    \text{PSL}(2,\mathbb{ R})$ in general.
Fortunately,
 the following theorem suggests that
  it does  not matter whether   $ B :\mathbb{R}/2  \mathbb{Z}\mapsto \text{SL}(2,\mathbb{R})$
or $ B :\mathbb{R}/  \mathbb{Z}\mapsto \text{PSL}(2,\mathbb{R})$ in defining reducibility.
 \end{remark}
 \begin{theorem}\label{XTh54}
 If $ B :\mathbb{R}/2  \mathbb{Z}\mapsto \text{SL}(2,\mathbb{R})$ is analytic in $|\Im x|<\eta$ and $ B(x+\alpha)^{-1}A(x)B(x)$
 is   constant, then  there exists  $ B' :\mathbb{R}/  \mathbb{Z}\mapsto \text{PSL}(2,\mathbb{R})$ being  analytic in $|\Im x|<\eta$  such that
 $ B'(x+\alpha)^{-1}A(x)B'(x)$ is  constant.
\end{theorem}
\textbf{Proof:}
\textbf{Step 1:} We will prove that there exists  $ B_1 :\mathbb{R}/4  \mathbb{Z}\mapsto \text{SL}(2,\mathbb{R})$ being  analytic in $|\Im x|<\eta$ such that
$ B_1^{-1}(x+\alpha)A(x)B_1(x)=V$ and  $ B_1(x+1)^{-1} B_1(x)=D$, where $ V$, $D$ are constant    and commute (i.e., $VD=DV$).
\par
By hypothesis  there exists $ B :\mathbb{R}/2  \mathbb{Z}\mapsto \text{SL}(2,\mathbb{R})$  such that  $ B^{-1}(x+\alpha)A(x)B(x)=V_1$ with $V_1$ being constant.
Let $D_1(x)=B(x+1)^{-1}B(x)$, then $D_1(x+2)=D_1(x)$ and
\begin{eqnarray}
\nonumber
   V_1D_1(x)V_1^{-1} &=&B(x+1+\alpha)^{-1}A(x+1)B(x+1) B(x+1)^{-1}B(x) B(x)^{-1}A(x)^{-1}B(x+\alpha) \\
   &=&  B(x+1+\alpha) ^{-1}B(x+\alpha)= D_1(x+\alpha). \label{XG5111}
\end{eqnarray}
 Assume that $V_1$ is not conjugate to a rotation  $R_\theta$ with
 $2\theta \in  \frac{1}{2}\alpha\mathbb{Z} +\mathbb{Z}$. Write $D_1(x)$ in the Fourier  series
  (note that $D_1(x)$ is well defined in $\mathbb{R}/2  \mathbb{Z}$)
\begin{equation}\label{XG511}
  D_1(x)=\sum_{k\in \mathbb{Z}}\hat{D }_1(k)e^{\pi i kx},\;\;\hat{D}_1(k)\in M(2,\mathbb{C}),
\end{equation}
then
\begin{equation}\label{XG512}
  \hat{D}_1(k)e^{\pi i k\alpha}=V_1\hat{D }_1(k)V_1^{-1}.
\end{equation}
If $\hat{D }_1(k)\neq 0   $ for some $k\neq 0$,  then $e^{\pi i k \alpha}$ is an eigenvalue of $\text{Ad}(V_1) : M(2, \mathbb{C})  \mapsto M(2,\mathbb{C})$, where $\text{Ad}(V_1)\cdot F=  V _1FV_1^{-1}$ for $F\in M(2, \mathbb{C})$.
 This implies that $V_1$ is conjugate to some rotation $R_\theta$  with $ 2\theta=\pm \frac{k\alpha}{2}+\ell $  for some $\ell \in \mathbb{Z}$
(see     Lemma \ref{XLe55} below), contradicting to our assumption. Thus   $D_1(x) = \hat{D }_1(0) $ is a constant.
Let $B_1(x)=B(x)$, $D=D_1$ and $V=V_1$,
then  $VD=DV$ by $(\ref{XG5111})$ and $ B_1^{-1}(x+\alpha)A(x)B_1(x)=V$.
\par
Assume that $ V_1$  is conjugate to some rotation  $R_\theta$  with $2\theta=\frac{k\alpha}{2}+\ell$,  where
$k, \ell \in \mathbb{Z} $,
i.e., $V_1=UR_\theta U^{-1}$ with $U\in \text{SL}(2,\mathbb{R}) $. Let $B_1(x)=B(x)UR_{\frac{k}{4}x}U^{-1}$, then $ B_1(x+4)=B_1(x)$
and
\begin{equation}\label{xx111111}
   B_1(x+\alpha)^{-1}A(x)B_1(x)=\pm I,
\end{equation}
where $I$ is the identity of $ 2\times2 $ matrix.
\par
  Let  $D_2(x)=B_1(x+1)^{-1}B_1(x) $. As $(\ref{XG5111})$, we have $D_2(x+\alpha)=D_2(x)$. By   the minimality of $ x\mapsto x+\alpha$,
$D_2$ is  constant. Let $V= \pm I$ and $D=D_2$, then    $VD=DV$   and $ B_1^{-1}(x+\alpha)A(x)B_1(x)=V$ by (\ref{xx111111}).
\par
\textbf{Step 2:} Let
$$d=\frac{1}{2\pi i}\int_{\Gamma}(\lambda I-\varepsilon D)^{-1}\ln\lambda d\lambda$$
 where $ \Gamma$ is a closed cure in complex plane,
contains all spectra of $\varepsilon D$ and $0\notin \Gamma$, and $
 \varepsilon \in \{-1,1\}$
($\varepsilon=1$ if the spectra of $D$ are positive,  otherwise $\varepsilon=-1$ ).
It is easy to check that
  $D=\varepsilon e^d$  and   $d\in \text{sl}(2,\mathbb{R})$  commutes with $V$ and $D$, where $d\in \text{sl}(2,\mathbb{R})$
  means the trace of matrix $d$ (denote $\text{tr} d $) is 0.
Let $B'(x)=B_1(x)e^{xd}$, then $B'(x+1)^{-1}B'(x)=\varepsilon \text{ I }$,   i.e.,  $ B' :\mathbb{R}/  \mathbb{Z}\mapsto \text{PSL}(2,\mathbb{R})$. Moreover,  $ B'(x+\alpha)^{-1}A(x)B'(x)=e^{-\alpha d}V$ is   constant.  $\qed$
\begin{lemma} \label{XLe55}
If for some  $k\in  \mathbb{Z} \backslash \{0 \}$ and $2\times 2$ matrix $D\neq 0$, the following holds,
\begin{equation} \label{XINjia}
   {D} e^{\pi i k\alpha}=V {D} V^{-1} ,
\end{equation}
where $V$ is a real constant cocycle.
Then  $V$ is conjugate to a rotation  $R_\theta$  with $ 2\theta=\pm \frac{k\alpha}{2}+\ell $ for some $\ell \in \mathbb{Z}$.
\end{lemma}
\textbf{Proof:} Without loss of generality, assume $ V $  is the form  of
$ \left(
    \begin{array}{cc}
       t & 0 \\
      0 & t^{-1} \\
    \end{array}
  \right)
$ with $ t\neq \pm1$ and $t \in \mathbb{R}$, or
$ \left(
    \begin{array}{cc}
      \pm 1 & a \\
      0 & \pm 1 \\
    \end{array}
  \right)
$ with $a\neq 0$ and $a\in \mathbb{R}$, or
$ \left(
    \begin{array}{cc}
      e^{2 \pi i\theta}   & 0 \\
      0 &   e^{-2 \pi i\theta}  \\
    \end{array}
  \right)
$ with  $\theta\in \mathbb{R}$, since $\det V=1$.
\par
If $ V =\left(
    \begin{array}{cc}
      e^{2 \pi i \theta}   & 0 \\
      0 &   e^{-2 \pi i\theta}  \\
    \end{array}
  \right) $.  Write  $D=(D_{ij})_{i,j=1,2}$, by a simple computation in (\ref{XINjia}),  we have

\begin{equation}\label{XINjia2}
   \left(
    \begin{array}{cc}
     D_{11}   &  D_{12}  \\
      D_{21}  &    D_{22}   \\
    \end{array}
  \right) e^{ \pi i k\alpha }=
  \left(
    \begin{array}{cc}
     D_{11}   & e^{4\pi i\theta} D_{12}  \\
     e^{-4\pi i\theta}D_{21}  &    D_{22}   \\
    \end{array}
  \right).
\end{equation}
Thus,   $D_{11}, D_{22}, D_{21}=0$ and $ e^{4\pi i\theta}=e^{ \pi i k\alpha }$; or
$D_{11}, D_{22}, D_{12}=0$ and $ e^{-4\pi i\theta}=e^{ \pi i k\alpha }$. In either case $ 2\theta=\pm \frac{k\alpha}{2}+\ell $ for
some $\ell \in \mathbb{Z}$.
\par
For
$V= \left(
    \begin{array}{cc}
       t & 0 \\
      0 & t^{-1} \\
    \end{array}
  \right)
$  with $t\neq\pm1$, or
$ V=\left(
    \begin{array}{cc}
      \pm 1 & a \\
      0 & \pm 1 \\
    \end{array}
  \right)
$ with $a\neq 0$,  we can prove that  those two cases can  not    happen   by a similar   discussion as the above. $\qed$

\begin{remark}\label{Re56}
By Theorem \ref{XTh54},  it does not matter whether    $ B :\mathbb{R}/2  \mathbb{Z}\mapsto \text{SL}(2,\mathbb{R})$
or $ B :\mathbb{R}/  \mathbb{Z}\mapsto \text{PSL}(2,\mathbb{R})$  in the   definition of reducibility.
The basic idea of the proof  in Theorem \ref{XTh54} is due to Avila-Krikorian $ \cite{AK}$,
where they
deal with another problem (Lemma 4.3, $ \cite{AK}$  ).
\end{remark}
\begin{lemma}\label{XLe417}
 Cocycle  $A=S_{\lambda,E}$ can not  be analytically  reducible to $\pm I$.
\end{lemma}
\textbf{Proof:} Otherwise,  without loss of generality,  assume there exists   $B:\mathbb{R}/\mathbb{Z}\rightarrow \text{PSL}(2,\mathbb{R})$
 being analytic, such that
 $B(x+\alpha)^{-1}A(x)B(x)=I
 $.
  Since $B(x)\in \text{PSL}(2,\mathbb{R})$,   $B(x+1)=\pm B(x)$. When $B(x+1)=B(x)$ the proof is simpler, see Remark \ref{XRe57}.
  Here we give the proof only for $B(x+1)=-B(x)$.
  Since  $B(x+\alpha)^{-1}A(x)B(x)=I$,
    it is easy to see that $ B$ must be with  the form $B(x)=\left(
                                                   \begin{array}{cc}
                                                     u_1(x) &  u_2(x) \\
                                                     u_1(x-\alpha)  &u_2(x-\alpha)  \\
                                                   \end{array}
                                                 \right)
  $
  and  \begin{eqnarray}
         (E-2\lambda\cos2\pi x) u_1(x)-u_1(x-\alpha) &=&  u_1(x+\alpha), \\
           (E-2\lambda\cos2\pi x) u_2(x) -u_2(x-\alpha) &=&  u_2(x+\alpha).
       \end{eqnarray}
      By comparing the Fourier coefficients (notice that both $u_1$ and $u_2$  are well defined in $\mathbb{R}/2\mathbb{Z} $),
      we obtain\begin{eqnarray}
                 (E-2\cos(\pi k\alpha)) \hat{u}_1(k)&=&   \lambda( \hat{u}_1(k+2)+ \hat{u}_1(k-2)), \label{X432}\\
                 (E-2\cos(\pi k\alpha)) \hat{u}_2(k)&=&   \lambda( \hat{u}_2(k+2)+ \hat{u}_2(k-2)),\label{X433}
               \end{eqnarray}
      where $\hat{u}_i(k)$ is the Fourier coefficients of $u_i$, $i=1,2$.
      \par
Let $ \tau $ be a new self-adjoint operator on $\ell^2(\mathbb{Z})$, with
\begin{equation}\label{X434}
  (\tau f) (k)=f(k+2)+f(k-2)+\frac{2}{\lambda}\cos(\pi k\alpha) f(k), \forall f\in \ell^2(\mathbb{Z}).
\end{equation}
After a simple computation
\begin{equation}\label{X435}
 \sum_{j=m}^{n}(f\tau g-g\tau f)(j)=W_n(f,g)-W_{m-1}(f,g),
\end{equation}
where
\begin{equation}\label{X436}
 W_n(f,g)=f(n)g(n+2) +f(n-1)g(n+1)-g(n)f(n+2)-g(n-1)f(n+1).
\end{equation}
In $(\ref{X435})$, let $f =\{\hat{u}_1(k) \}_{k\in \mathbb{Z}} $ and $g=\{\hat{u}_2(k) \}_{k\in \mathbb{Z}}$, combining with $  (\ref{X432}) $ and $ ( \ref{X433}) $,
one has
\begin{equation}\label{X437}
  W_n(\hat{u}_1,\hat{u}_2)=W_{m }(\hat{u}_1,\hat{u}_2).
\end{equation}
Since  $ u_i$ is analytic,      $\lim_{n\rightarrow \infty}\hat{u}_i(n)=0$ for   $i=1,2 $ and $\lim_{m\rightarrow \infty} W_m(\hat{u}_1,\hat{u}_2)=0$.
By (\ref{X437}),
\begin{equation}\label{X438}
 W_n(\hat{u}_1,\hat{u}_2)=\hat{u}_1(n)\hat{u}_2(n+2) +\hat{u}_1(n-1)\hat{u}_2(n+1)-\hat{u}_2(n)\hat{u}_1(n+2)-\hat{u}_2(n-1)\hat{u}_1(n+1)=0.
\end{equation}
  Moreover,  $\hat{u}_i(k)=0$ for even $k$ because of $u_i(x+1)=-u_i(x)$, $i=1,2 $.
 In $(\ref{X438})$, let $n=2k$,
we have
\begin{equation}\label{X439}
\hat{u}_1(2k-1)\hat{u}_2(2k+1)-\hat{u}_2(2k-1)\hat{u}_1(2k+1)=0.
\end{equation}
This implies $\hat{u}_1$ and $\hat{u}_2$ are linear related,  contradicting to $\det B=1$. $\qed$
\begin{remark}\label{XRe57}
For  another  case $B(x+1)=B(x)$, i.e.,    $B:\mathbb{R}/\mathbb{Z}\mapsto  \text{SL}(2,\mathbb{R})$, the proof is simper. We only need replace (\ref{X434})
with $(\tau f) (k)=f(k+1)+f(k-1)+\frac{2}{\lambda}\cos(2\pi k\alpha) f(k)$  and (\ref{X436}) with
$W_n(f,g)=f(n)g(n+1)-g(n )f(n+1)$.
\end{remark}

If $E\in \Sigma _{\lambda,\alpha}$ such that $\theta(E)$    is not   $\epsilon_0$-resonant,
by Remark $ \ref{xxRe315}$, there
  exists a non-zero exponentially decaying solution   of  $\hat{H}\hat{u}=E\hat{u}$.  Next  we will set up the
  reducibility of cocycle $A=S_{\lambda,E}$  via  constructing  reducible matrix.
  \begin{theorem}\label{Th417}
  Given $\alpha\in \mathbb{R}\backslash \mathbb{Q }$, $\theta\in \mathbb{R}$ and $E\in  \Sigma_{\lambda,\alpha}$, suppose  there
  exists a non-zero exponentially decaying eigenfunction $ \hat{u}=\{\hat{u}_{k}\}_{k\in \mathbb{Z}}$, i.e.,   $\hat{H}_{\lambda,\alpha,\theta}\hat{u}=E\hat{u}$ with  $ |\hat{u}_k|\leq   e^{-2\pi \eta|k|} $ for $k$ large enough, then
  the following hold.
  \par
  (1)If $2\theta\notin \alpha \mathbb{Z }+ \mathbb{Z }$, then  there exists $B:\mathbb{R}/\mathbb{Z}\mapsto \text{SL}(2,\mathbb{R})$
  being analytic
  in $|\Im x|<  \eta $, such that $B(x+\alpha)^{-1}A(x) B(x)=R_{\pm\theta}$ , i.e., $(\alpha,A)$   is analytically reducible in trip  $ |\Im x|<\eta$,     where $A=S_{\lambda,E}$. In
  this case $ \rho(\alpha,A)=\pm\theta+\frac{m}{2}\alpha \mod \mathbb{Z}$, where $m=deg(B)$\footnote{Since $B$ is well defined in $\mathbb{R}/\mathbb{Z}$,
      $m=$deg($B$) must be even  . } and $2 \rho(\alpha,A)\notin \alpha \mathbb{Z }+ \mathbb{Z }$.
  \par
  (2)If $2\theta-k\alpha\in\mathbb{ Z}$ for some $k\in\mathbb{ Z}$ and   $\eta> 8\beta(\alpha)$, then there exists
   $B:\mathbb{R}/\mathbb{Z}\mapsto \text{PSL}(2,\mathbb{R})$  being analytic
  in $|\Im x|< \frac{\eta}{4}$, such that $B(x+\alpha)^{-1}A(x) B(x)= \left(
                                                         \begin{array}{cc}
                                                           \pm1 & a \\
                                                           0 & \pm 1\\
                                                         \end{array}
                                                       \right)
  $  with $a\neq 0$,    i.e., $(\alpha,A)$   is analytically reducible in trip   $|\Im x|<\frac{\eta}{4}$.  In
  this case $2\rho(\alpha,A)=   m\alpha \mod  \mathbb{Z}$, where $m=\deg(B)$, i.e., $2 \rho(\alpha,A)\in \alpha \mathbb{Z }+ \mathbb{Z }$.
  \end{theorem}
 \textbf{ Proof:}
 Since $ |\hat{u}_k|\leq   e^{-2\pi \eta|k|} $ for $k$ large enough, $u(x)=\sum \hat u_k e^{2 \pi i k x}$ is analytic
  in $|\Im x|< \eta$.
 Let
   $U(x)= \left(
           \begin{array}{c }
             e^{2 \pi i \theta} u(x)\\
            u(x-\alpha) \\
           \end{array}
         \right)
  $, then (see \S 2.4)
  \begin{equation}\label{xxxx222}
   A(x) \cdot U(x)=e^{2
\pi i \theta} U(x+\alpha).
\end{equation}
Let $\tilde B(x)$ be a matrix with columns
$U(x)$ and $\overline U(x)$, i.e., $ \tilde B(x)=\left( U(x), \overline U(x)\right)$. Note that $\overline U(x)$ is given by  footnote 6.
By the minimality of $x \mapsto x+\alpha$ and
(\ref{xxxx222}),
$\det \tilde B$ is a constant.
\par
(Case A)
If $\det \tilde B \neq 0$, we have  $\tilde B(x+\alpha)^{-1} A(x) \tilde
B(x)=\left(
       \begin{array}{cc}
         e^{2 \pi i \theta} &0 \\
         0 & e^{-2 \pi i \theta} \\
       \end{array}
     \right)
$.
It is easy to see that $\det \tilde B=\pm c i$ for some $c>0$. If  we  take
$B =\frac{1}{(2c)^{1/2}}\tilde{B}\left(
                                       \begin{array}{cc}
                                         1 & \pm i \\
                                         1 & \mp i \\
                                       \end{array}
                                     \right)
$,  then $B(x+\alpha)^{-1}A(x) B(x)=R_{\pm\theta}$, and  $\rho(\alpha,A)=\pm\theta+\frac{m}{2}\alpha \mod \mathbb{Z}$, where $m=\deg(B)$.
 \par
Now we are in position to prove  $2\theta\notin \alpha \mathbb{Z }+ \mathbb{Z }$. Otherwise, there exists some $k\in \mathbb{Z}$ such that  $2\theta-k\alpha\in \mathbb{Z}$.
 Let $B'(x)=B(x)R_{\pm \frac{kx}{2}}$, we have $B'(x+\alpha)^{-1}A(x) B'(x)= \text{I}$ or $B'(x+\alpha)^{-1}A(x) B'(x)= -\text{I}$ . This is impossible by Lemma $\ref{XLe417}$.
 \par
 (Case B) If $\det \tilde B=0$.
 By the
minimality of $x \mapsto x+\alpha$ and (\ref{xxxx222}), $U(x) \neq 0$ for all $x$ with $ |\Im x|<\eta$.  Thus we have
$U(x)=\psi(x) W(x)$ with $W(x+1) =\pm W(x)$  and $W(x) $  being    real analytic
  in $|\Im x|< \eta$,   and $|\psi(x)|=1$ for $x\in\mathbb{ R}$ (see    Lemma \ref{xxLe510} below).

  There exists $\delta>0$ such that  $ ||W(x)||>\delta $ in  $|\Im x|<\frac{\eta}{2}$, since $W(x)\neq 0$ for all $x$ with $|\Im x|<\eta$. Let $B_1$  be given by Lemma
  $\ref{matrix B}$ with first column $W$, then $ B_1 :\mathbb{R}/ 2\mathbb{Z}\mapsto \text{SL}(2,\mathbb{R})$  is analytic  in $|\Im x|< \frac{\eta}{2}$,  and $B_1(x+\alpha)^{-1}A(x)B_1(x)=
\left(
  \begin{array}{cc}
    d(x) & \kappa(x)\\
     0& d(x)^{-1} \\
  \end{array}
\right)
 $ with  $d(x)=\frac{\psi(x+\alpha)}{\psi(x)}e^{2\pi i\theta}$.
 Since
 $|d(x)|=1$ and $d(x)$ is real for $x\in \mathbb{R}$,   $d(x)=\pm 1$, i.e.,
 $B_1(x+\alpha)^{-1}A(x) B_1(x)= \left(
                                                         \begin{array}{cc}
                                                           \pm1 & \kappa(x) \\
                                                           0 & \pm 1\\
                                                         \end{array}
                                                       \right)$.
Moreover,  $2\rho(\alpha,A)= m_1\alpha \mod \mathbb{Z}$ since the degree of $ \left(\begin{array}{cc}
                                                           \pm1 & \kappa(x) \\
                                                           0 & \pm 1\\
                                                         \end{array}
                                                       \right) $ is 0, where  $m_1=\deg(B_1)$.
\par
If $ \eta> 8\beta $,
we can   further conjugate $A$ to a constant parabolic
matrix by solving (comparing Fourier coefficients) the  homological equation
$\pm \phi(x+\alpha) \mp \phi(x)=\kappa(x)-\int_0^2
\kappa(x) dx$ in $ \mathbb{R}/ 2\mathbb{Z}$ with $  \hat{\phi}_0=0$. More precisely, $ \hat{\phi}_k=\mp\frac{\hat{\kappa}_k}{1-e^{ \pi i k\alpha}}$, $k\neq 0$, thus $\phi$
is analytic in $|\Im x|<\frac{\eta}{4}$ because of $\kappa(x)$ being  analytic in $|\Im x|<\frac{\eta}{2}$ and small divisor condition (\ref{X31}).
  Let  $B_2(x)= B_1(x)
\left(
  \begin{array}{cc}
    1 &  \phi(x) \\
    0 & 1 \\
  \end{array}
\right)
 $, we get
$B_2(x+\alpha)^{-1} A(x) B_2(x)=
\left(
  \begin{array}{cc}
    \pm 1 & \int_0^2 \kappa(x) dx \\
    0 & \pm 1\\
  \end{array}
\right)
$  and  $ B_2  $ is well defined in  $\mathbb{R}/ 2\mathbb{Z}$.
\par
By Theorem $\ref{XTh54}$ (let $B=B_2$ in  Theorem $\ref{XTh54}$),  there exists   $ B_3 :\mathbb{R}/  \mathbb{Z}\mapsto \text{PSL}(2,\mathbb{R})$ such  that
$B_3(x+\alpha)^{-1} A(x) B_3(x)$ is a constant cocycle $C$. We will prove that $C$ is conjugate to
$\left(
  \begin{array}{cc}
    \pm 1 & a   \\
    0 & \pm  1\\
  \end{array}
\right)$ with $a$ a constant. Otherwise, $C$ is conjugate to rotation $R_{\theta'}$ with $ 2\theta' \in  \alpha\mathbb{Z}+\mathbb{Z}$ (since $2\rho(\alpha,A)\in   \alpha \mathbb{Z}  + \mathbb{Z}$), this is impossible by the discussion in Case A; or  $C$ is conjugate to
$\left(
  \begin{array}{cc}
    t & 0   \\
    0 &  t^{-1}\\
  \end{array}
\right)$ with $t\neq \pm 1$, this is impossible since $E\in \Sigma_{\lambda,\alpha}$ ($S_{\lambda, E}$ is not uniformly hyperbolic\footnote{ We say that
cocycle $(\alpha,A) $ is uniformly hyperbolic if there exist constants $c > 0$,
$\gamma>1$ such that $||A_n(x)||\geq c\gamma^n$ for every $x\in \mathbb{R}$  and $n >0$.} for $E\in \Sigma_{\lambda,\alpha}$, see $ \cite{Jo}$).
 Therefore, there exists a cocycle $U$ such that $U^{-1}CU= \left(
  \begin{array}{cc}
    \pm 1 & a   \\
    0 & \pm  1\\
  \end{array}
\right)  $. Let $B(x)= B_3(x)U$, then  $B(x+\alpha)^{-1}A(x) B(x)= \left(
                                                         \begin{array}{cc}
                                                           \pm1 & a \\
                                                           0 & \pm 1\\
                                                         \end{array}
                                                       \right)$.
                                                       This implies that $2\rho(\alpha,A)= m \alpha \mod \mathbb{Z}$, where  $m =\deg(B)$.
                                                       Notice that   $a\neq 0$  by Lemma $\ref{XLe417}$.
\par
Now we prove that $2
\theta=k\alpha \mod \mathbb{Z}$.
Since $d=\pm 1$ and $d(x)=\frac{\psi(x+\alpha)}{\psi(x)}e^{2\pi i\theta}$,
 $ \psi(x+\alpha)e^{2 \pi i \theta}=\pm   \psi(x) $.
This implies ( comparing Fourier coefficients) that
$\psi(x)=e^{-\pi i k x}$ (notice that $\psi$ is well defined   in
$ \mathbb{R}/2 \mathbb{Z}$) and $e^{2 \pi i \theta}=\pm e^{\pi i k \alpha}$ for some $k\in \mathbb{Z}$, that is $2
\theta=k\alpha \mod \mathbb{Z}$.
 \par
 Putting case A and   B together, we finish the proof.
 \begin{remark}\label{Re57}
 In above discussion, we have prove: if $(\alpha,A)$ is reducible  and $2\rho(\alpha,A)\in \alpha\mathbb{ Z}+\mathbb{Z}$,
 where $ A=S_{\lambda,E}$  with $E\in \Sigma_{\lambda,\alpha}$, then $(\alpha,A)$ must be  conjugate to
$\left(
  \begin{array}{cc}
    \pm 1 & a  \\
    0 & \pm 1\\
  \end{array}
\right)$ with $a\neq 0$.
 \end{remark}

\begin{lemma}\label{xxLe510}
Under the notation of Theorem  \ref{Th417},
if $\det \tilde B=0$,  we have
$U(x)=\psi(x) W(x)$ with $W(x) $ being    real analytic
  in $|\Im x|< \eta$ and $W(x+1) =\pm W(x)$,   and $|\psi(x)|=1$ for $x\in \mathbb{R}$.
\end{lemma}
\textbf{Proof:}
Let $U(z)= \left (\begin{array}{c}
             u_1(z) \\
              u_2(z)
           \end{array}\right )
$, $|\Im z|< \eta$. By condition $\det \tilde B=0$, then there exists $k(z)$ such that
\begin{equation}\label{G430}
   u_1(z)=k(z)\overline{u}_1(z) \text{ and }
 u_2(z)=k(z)\overline{u}_2(z) .
\end{equation}
  By
minimality of $z \mapsto z+\alpha$ and (\ref{xxxx222}), $U(z) \neq 0$ for   $ |\Im z|<\eta$. Thus  $k(z)\neq 0$ for all $ |\Im z|<\eta$.
Moreover, $k(z)=\frac{u_1(z)}{\overline{u}_1(z)}$ or $k(z)=\frac{u_2(z)}{\overline{u}_2(z)}$, which implies
$k(z)$  can be selected so that  $k(z)$  is   analytic in $ |\Im z|<\eta$ and $|k(x)|=1$ for $x\in \mathbb{R}$.
\par
We will prove that there exists $\varphi$ being   analytic     in $ |\Im z|<\eta$ such that $\varphi^2=k$ and $\varphi$ is well
defined in $\mathbb{R}/2\mathbb{Z}$  with $|\varphi(x)|=1$ for $x\in \mathbb{R}$ (i.e., $ \overline{\varphi}\varphi=1$).  Fix a point  $z_0\in \mathbb{R}$,   and solve  $p'(z)=\frac{k'(z)}{k(z)}$ with  $p(z_0)=\ln k(z_0)$ ( selecting  a branch).
We have $p(z)$ is  analytic in $ |\Im z|<\eta$ and $e^{p(z)}=k(z)$. Let $\varphi(z)=e^{\frac{1}{2}p(z)}$, then $\varphi^2=k$. By the uniqueness theorem of analytic function in Complex Analysis, it's easy to verify  $\varphi$ is well
defined in $\mathbb{R}/2\mathbb{Z}$  and $|\varphi(x)|=1$ for $x\in \mathbb{R}$. Combining  with $ (\ref{G430})$, for  $x\in \mathbb{R}$,
$\overline{\varphi}(x)u_i(x)=\varphi(x)\overline{u_i}(x)$, $ i=1,2$,
which implies both $\overline{\varphi}u_1$ and $\overline{\varphi}u_2$ are real analytic in $ |\Im z|<\eta$.
Let $W=\left(\begin{array}{c}
          \overline{\varphi}u_1 \\
         \overline{\varphi}u_2
       \end{array}
       \right)
$ and $ \psi=\varphi$, we prove this lemma.
  \begin{theorem} \label{Th420}
   For $0<\beta(\alpha)<\infty$ and $|\lambda|<e^{-C_2\beta}$,  let $A=S_{\lambda,E}$ with  $E\in \Sigma_{\lambda,\alpha}$.
   If $2\rho(\alpha,A)\in \alpha\mathbb{Z} +\mathbb{Z}$ (i.e., $ N_{\lambda,\alpha} (E)\in \alpha\mathbb{Z} +\mathbb{Z}$ by (\ref{G26})),   then there exists
  $B : \mathbb{R}/\mathbb{Z } \mapsto \text{PSL}(2,\mathbb{R})$,
analytically extending to $ |\Im x| < \frac{h}{4}$, such that
$B(x+\alpha)^{-1}A(x) B(x)= \left(
                                                         \begin{array}{cc}
                                                           \pm1 & a \\
                                                           0 & \pm 1\\
                                                         \end{array}
                                                       \right)
  $,
  with $a\neq 0$.
  \end{theorem}
\textbf{Proof:} Let $E\in \Sigma_{\lambda,\alpha}$, we first prove that if $2\rho(\alpha,A)\in \alpha\mathbb{Z} +\mathbb{Z}$, then
  $\theta(E)$ given by Lemma  \ref{Le41} is  not $\epsilon_0$-resonant. Otherwise,
  by Theorem $\ref{Th415}$, there exists
   $m_j$ such that $|m_j|<C|n_j|$ and $ ||2 \rho(\alpha,A) -m_j\alpha\pm(2  \theta -n_j\alpha ) ||_{\mathbb{R}/ \mathbb{Z }}<e^{-ch |n_{j+1}|}$.
    By $(\ref{XG42})$
    \begin{equation}\label{G431}
      ||2 \rho(\alpha,A) -m_j\alpha||_{\mathbb{R}/ \mathbb{Z }}\geq ||  2  \theta -n_j\alpha  ||_{\mathbb{R}/ \mathbb{Z }}- e^{-ch |n_{j+1}|} >e^{-8\beta |n_{j+1}|}-e^{-ch |n_{j+1}|}>0,
    \end{equation}
and
  \begin{equation}\label{G432}
  || 2\rho(\alpha,A) -m_j\alpha||_{\mathbb{R}/\mathbb{Z }}\leq ||  2  \theta -n_j\alpha  ||_{\mathbb{R}/\mathbb{Z }}+ e^{-ch |n_{j+1}|}
 \leq e^{- \epsilon_0|n_j|}+e^{-ch |n_{j+1}|}\leq e^{-c\epsilon_0|m_j|}.
  \end{equation}
It follows from (\ref{G432}) that  $\rho(\alpha,A)$ has a $c\epsilon_0$-resonance at $m_j$ if  $|m_j|$ is large enough by Lemma \ref{XXXLe33}.
 If the set of $c\epsilon_0$-resonance for $\rho(\alpha,A)$ is finite, i.e., $\{m_j\}$ is finite,
 by $(\ref{G431})$,  there exists some $\delta>0$ such that $|| 2\rho(\alpha,A) -m_j\alpha||>\delta$  for all $j$, which is  contradicted to the
 fact $ || \rho(\alpha,A) -m_j\alpha||_{\mathbb{R}/\mathbb{Z }} \rightarrow0 \text{ as } {j\rightarrow \infty}$ by  the second  inequality
 in (\ref{G432}). Thus $\rho(\alpha,A)$ is   $c\epsilon_0$-resonant,  this is impossible because
 of $2\rho(\alpha,A)\in \alpha\mathbb{Z} +\mathbb{Z}$.
 \par
  Now that  $\theta(E)$ is not $\epsilon_0$-resonant, by Remark $\ref{xxRe315}$
there
  exists a non-zero exponentially decaying solution $\hat u$ of $  \hat{H}_{\lambda,\alpha,\theta}\hat{u}=E\hat{u}$ with   $|u_k|\leq e^{-2\pi h|k|} $ for $|k|$ large enough,
  where $h=C_1^2 \beta$ by our hypothesis in the beginning of \S 4.
  Combining with Theorem  $\ref{Th417}$,
 we finish the proof.
  \par

   \begin{theorem} \label{Th421}$(\text{Theorem 4.1}, \cite{AJ1})$
  For $ \beta(\alpha)=0 $ and $|\lambda|<1$,  let $A=S_{\lambda,E}$ with $E\in \Sigma_{\lambda,\alpha}$,  there exists a small
  constant $ c(\lambda,\alpha)$ such that,  if $2\rho(\alpha,A)\in \alpha \mathbb{Z} + \mathbb{Z} $ then there exists $B : \mathbb{R}/\mathbb{Z } \mapsto \text{PSL}(2,\mathbb{R})$ being
  analytic in $ |\Im x| <c(\lambda,\alpha)$ such that $B(x + \alpha )^{-1}A(x)B(x)$ is   constant.

  \end{theorem}

\begin{remark}
  Avila and Jitomirskaya prove Theorem $ \ref{Th421}$ only for $\alpha \in DC$, in fact, their proof suggests it holds for all
  $\beta(\alpha)=0$ (after carefully checking their proof).
\end{remark}
\section{Proof of Main Theorem}
\begin{theorem}\label{Th51}
  If $E_0\in \Sigma_{\lambda,\alpha}$  such that $  2\rho(\alpha,A_{E_0})\in \alpha \mathbb{Z} + \mathbb{Z}$, and $(\alpha,A_{E_0})$ is analytically
  reducible in ${|\Im x| < \eta}$  with $\eta>6\beta(\alpha)$  $(0\leq \beta(\alpha)<\infty)$, where $A_{E_0}=S_{\lambda,E_0}$, then $E_0$ is an     endpoint
of some   gap.
  \end{theorem}
\textbf{Proof:}  Here we only give the proof  if  $0< \beta(\alpha)<\infty$.  For $\alpha$ with $ \beta(\alpha)=0$, the proof is similar.
  Let $B : \mathbb{R}/\mathbb{Z } \mapsto \text{PSL}(2,\mathbb{R})$ be
  analytic in ${|\Im x| < \eta}$ such that $B(x + \alpha )^{-1}A_{E_0}(x)B(x)$ is a constant cocycle. Since $  2\rho(\alpha,A_{E_0})\in \alpha \mathbb{Z} + \mathbb{Z}$,  combining with Remark \ref{Re57}, we have
     \begin{equation}\label{G51}
     B(x + \alpha )^{-1}A_{E_0}(x)B(x)=
         \left(
           \begin{array}{cc}
           \pm 1& a \\
             0 & \pm1  \\
           \end{array}
         \right),
     \end{equation}
    with $a \neq 0$. Without loss of generality,   assume
    $ B(x + \alpha )^{-1}A_{E_0}(x)B(x)=
         \left(
           \begin{array}{cc}
            1& a \\
             0 &  1  \\
           \end{array}
         \right)\triangleq  Z$ with  $a<0$. Writing $B=(B_{ij})_{i,j=1,2}$,   one easily obtains the following facts,
   \begin{equation}\label{G52}
     B_{21}(x+\alpha)=B_{11}(x),  B_{22}(x+\alpha)=B_{12}(x)-aB_{21}(x+\alpha) .
   \end{equation}
   Below, let $\varepsilon>0$ be small.  After carefully computing,
     \begin{equation}\label{G53}
     B(x + \alpha )^{-1}A_{E_0+\varepsilon}(x)B(x)=Z+\varepsilon P,
     \end{equation}
      where
     \begin{equation}\label{G54}
     P=
        \left(
          \begin{array}{cc}
            B_{11}B_{12}-aB^2_{11} & -aB_{11}B_{12}+B_{12}^2\\
            -B_{11}^2 &  - B_{11}B_{12} \\
          \end{array}
        \right).
     \end{equation}
      We will prove that for an appropriate cocycle $B_1   :\mathbb{R}/\mathbb{Z}\mapsto \text{SL}(2,\mathbb{R})  $, one has
     \begin{equation}\label{G55}
     B_1(x + \alpha )^{-1}(Z +\varepsilon P(x))B_1(x)=Z+\varepsilon [P]+O(\varepsilon^2),
     \end{equation}
     where $[\cdot ]$ denotes the average of a  matrix-valued function  over $\mathbb{R}/\mathbb{Z}$.
     This can be done by a step KAM iteration (or averaging theory). Refer to \cite{Amo}. Namely,  we will look for a cocycle  $B_1$ with the form of
        $ B_1= e^{\varepsilon Y} $, where $Y   :\mathbb{R}/\mathbb{Z}\mapsto \text{sl}(2,\mathbb{R})  $   (i.e., $Y(x+1)=Y(x)$ and $ tr (Y(x))=0$).
        Clearly,
     \begin{eqnarray}
     \nonumber
       B_1(x + \alpha )^{-1}(Z +\varepsilon P(x))B_1(x) &=&  (I-\varepsilon Y(x  + \alpha )+O(\varepsilon^2)) (Z+\varepsilon P )(I+\varepsilon Y+O(\varepsilon^2)) \\
         &=&  Z+\varepsilon(ZY(x)+P(x)-Y(x+\alpha) Z )+O(\varepsilon^2) .\label{G56}
     \end{eqnarray}
  Let $T(x)=Z^{-1} P(x)-\frac{tr(Z^{-1} P)}{2}I$ and
  solve the  homological  equation
  \begin{equation}\label{1}
    Y(x+\alpha)Z-ZY(x)=Z(T(x)-\hat{T}(0)) \text{ in } \mathbb{R}/\mathbb{Z}
  \end{equation}
  with $\hat{Y}(0)=0$.
   We get
   $\hat{Y}_{11}(k) = \frac{\ast}{(1-e^{2\pi ik \alpha})^2}$, $\hat{Y}_{12}(k) = \frac{\ast}{(1-e^{2\pi ik \alpha})^3}$,
  $\hat{Y}_{ 2 1}(k) = \frac{\ast}{(1-e^{2\pi ik \alpha})}$ and $\hat{Y}_{22}(k) = \frac{\ast}{(1-e^{2\pi ik \alpha})^2}$, $ k\neq 0$,
  where  $\hat{Y}_{ij}(k)$ is the Fourier  coefficients of matrix elements $Y_{ij}$ of $Y$, $i,j=1,2$,  and  $\ast$  may be different. Using small divisor condition (\ref{X31}),     $Y$ is  analytic  if $\eta>6\beta$.
  Since $Y$ is a solution of equation $Y (x +\alpha)-ZY(x)Z^{-1} = Z(T  (x)- \hat{T} (0))Z^{-1}$,
$tr(  Y (x +\alpha))-tr (Y(x)) = tr(T  (x)-  \hat{T} (0)) = 0$, i.e.,  $tr Y(x) $  is constant for   $ x  \in  \mathbb{R}/ \mathbb{Z}$.
Notice that $ \hat{Y} (0)  =  0$, then $tr(Y(x))=0$  for $ x  \in  \mathbb{R}/ \mathbb{Z}$, i.e., $B_1=e^{\varepsilon Y}$ is indeed a cocycle.
\par
  By (\ref{G53})  $\det(Z+\varepsilon P)=1$, it is  straightforward to compute that    $tr(Z^{-1} P)=-\varepsilon \det P $,
  thus the coefficients of $\varepsilon$ in  $ (\ref{G56})$ satisfies
    \begin{equation}\label{1}
    ZY(x)+P(x)-Y(x+\alpha) Z=[P]+O(\varepsilon),
    \end{equation}
which implies $(\ref{G55})$.
\par
  Moreover,
\begin{equation}\label{G59}
   Z+\varepsilon[P]+O(\varepsilon^2)=\exp(Z_0+\varepsilon Z_1+O(\varepsilon ^2)) ,
\end{equation}
where
\begin{equation}\label{G60}
Z_0=
     \left(
       \begin{array}{cc}
         0 & a \\
         0 & 0 \\
       \end{array}
     \right), \text { and }
     Z_1=
     \left(
          \begin{array}{cc}
            [B_{11}B_{12}]-\frac{a}{2}[B^2_{11}] & -a[B_{11}B_{12}]+[B_{12}^2]\\
            -[B_{11}^2] &  - [B_{11}B_{12}]+\frac{a}{2}[B^2_{11}] \\
          \end{array}
        \right).
\end{equation}
Let
\begin{equation*}
D=
     \left(
       \begin{array}{cc}
         d_1 & d_2 \\
         d_3 & -d_1 \\
       \end{array}
     \right)=Z_0+\varepsilon Z_1,
\end{equation*}
whose determinant is $d=-d_1^2-d_2d_3=a\varepsilon [B_{11}^2]+O(\varepsilon^2)<0$  for  small $ \varepsilon>0$, since $ [B_{11}^2]>0$
(otherwise $B_{11}=0$, by $ (\ref{G52})$ $B_{21}=0$, this is impossible ). Now we
let
\begin{equation*}
F=
     \left(
       \begin{array}{cc}
         d_2 & d_2 \\
         -d_1+ \sqrt {-d} & -d_1- \sqrt {-d} \\
       \end{array}
     \right),
\end{equation*}
 which has determinant
$-2a\sqrt {-a\varepsilon [B_{11}^2]}+O(\varepsilon) $, then $||F||=O(1)$, $||F^{-1}||=O(\varepsilon^{-1/2})$,  and
\begin{equation*}
     F^{-1}DF=\left(
       \begin{array}{cc}
         \sqrt[ ]{-d} &0 \\
        0 &  -\sqrt[ ]{-d} \\
       \end{array}
     \right)\triangleq H.
\end{equation*}
Moreover,
\begin{equation}\label{1222}
    \exp(Z_0+\varepsilon Z_1+O(\varepsilon^2))=\exp(F(H+O(\varepsilon^{3/2}))F^{-1})=F\exp( H+O(\varepsilon^{3/2}) )F^{-1}.
\end{equation}
Notice  that
\begin{equation*}
    H+O(\varepsilon^{3/2})= \sqrt[ ]{-a\varepsilon  [B_{11}^2]} \left(
    \left(
      \begin{array}{cc}
        1 & 0 \\
        0 & -1\\
      \end{array}
    \right)
    +O(\varepsilon )\right).
\end{equation*}
Therefore,
if $ \varepsilon$ is small enough, the cocycle $ A_{E_0+\varepsilon}$  has an exponential dichotomy
( i.e.,  $ A_{E_0+\varepsilon}$ is   uniformly hyperbolic ),  which implies
$E_0+\varepsilon \notin \Sigma_{\lambda,\alpha}$, i.e., $E_0$ is an endpoint of some   gap. $\qed$

  \begin{remark}
  In   $\cite{P2}$,  Puig proves Theorem $\ref{Th51}$ for $\alpha\in DC $,  we extend   his result to   all $\alpha$ with $\beta(\alpha)<\infty $.
  \end{remark}
  Combining with Avila and Jitomirskaya's work \cite{AJ1},\cite{AJ2}, we give a summary  of the dry version of Ten Martini Problem.
\begin{theorem}\label{Conclution}
For every $ \alpha\in \mathbb{R}\backslash \mathbb{Q}$, let $\beta(\alpha)$ be given by (\ref{G11}), then the following statements  hold.
\par
(1) If $\beta(\alpha)=\infty$, then $\Sigma_{\lambda,\alpha}$ has all gaps open for all $\lambda  \neq 0$.
\par
(2) If $0<\beta(\alpha)<\infty$, then $\Sigma_{\lambda,\alpha}$ has all gaps open for   $ 0<|\lambda|< e^{-C_2\beta}$, or $e^{-\beta}<|\lambda|<e^{\beta}$,
  or $  |\lambda|> e^{ C_2\beta}$, where $C_2$ is a large absolute  constant.
\par
(3)If  $\beta(\alpha)=0$, then $\Sigma_{\lambda,\alpha}$ has all gaps open if $\lambda  \neq 0,-1,1$.
\end{theorem}
  \textbf{Proof :}
  If $ \beta(\alpha)=\infty$,  this case has already been  proved by Avila and Jitomirskaya  (Theorem 8.2, $\cite{AJ1}$).
  \par
  If $0 <\beta(\alpha)<\infty$,  Avila and Jitomirskaya  (Theorem 8.2, $\cite{AJ1}$) have proved  that  $\Sigma_{\lambda,\alpha}$ has all gaps open for $e^{-\beta}<|\lambda|<e^{\beta}$.
   Fix $ \epsilon_0=C_1\beta$, $h=C_1\epsilon_0$, where $C_1$ is a large absolute constant given in Theorem \ref{Th32}.
    Let $C_2$ be a large absolute constant  also given  in the beginning of \S 4.    If  $ |\lambda|<e^{-C_2\beta}$,  by Theorem $\ref{Th420}$, for any spectrum
   $   E_0$ satisfying  $  N_{\lambda,\alpha}(E_0)\in \alpha \mathbb{Z} + \mathbb{Z}$, i.e., $  2\rho(\alpha,A_{E_0})\in \alpha \mathbb{Z} + \mathbb{Z}$, there exists $B:\mathbb{R}/\mathbb{Z } \mapsto  \text{PSL}(2,\mathbb{R})$ being
  analytic in ${|\Im x| < \frac{h}{4}}$ such that $B(x + \alpha )^{-1}A_{E_0}(x)B(x)$ is   constant. Notice that  $ \frac{h}{4}>6\beta$,  since $C_1$ is large.
  By  Theorem $\ref{Th51} $,  $E_0$ is an     endpoint
of some   gap.
   For $ |\lambda|> e^{ C_2\beta}$, notice that
  $\Sigma_{\lambda^{-1},\alpha}=\lambda^{-1}\Sigma_{\lambda,\alpha}$ and
  $N_{\lambda^{-1},\alpha}(\lambda^{-1}E)= N_{\lambda,\alpha}(E)$ (Aubry duality).
  \par
  If $\beta(\alpha)=0$, we only need  replace      Theorem $\ref{Th420}$ with Theorem $\ref{Th421}$.

           \begin{center}
           
             \end{center}
  \end{document}